\documentclass[11pt]{article}
\usepackage[a4paper,margin=1in]{geometry}
\usepackage{amsmath,amssymb,amsthm,mathtools}
\usepackage{bm}
\usepackage{graphicx}
\usepackage{booktabs}
\usepackage{multirow}
\usepackage{xcolor}
\usepackage{hyperref}
\usepackage{enumitem}
\usepackage{algorithm}
\usepackage{algorithmic}
\usepackage{listings}
\usepackage{tabularx}
\usepackage{subcaption}
\usepackage{MnSymbol}
\usepackage{tikz}
\usetikzlibrary{arrows.meta,positioning,calc}
\usepackage{authblk}

\hypersetup{
  colorlinks=true,
  linkcolor=blue,
  citecolor=blue,
  urlcolor=blue
}

\DeclareMathOperator{\herm}{herm}
\newcommand{\HH}{\mathbb{H}}
\newcommand{\R}{\mathbb{R}}

\newcommand{\Herm}{\mathrm{Herm}}
\newcommand{\ip}[2]{\left\langle #1,#2\right\rangle_{\R}}
\newcommand{\tr}{\mathrm{tr}}
\newcommand{\Rea}{\operatorname{Re}}

\lstdefinestyle{quatica}{
  basicstyle=\ttfamily\small,
  keywordstyle=\color{blue},
  commentstyle=\color{gray},
  stringstyle=\color{teal},
  showstringspaces=false,
  breaklines=true,
  frame=single,
  rulecolor=\color{black!20},
  columns=fullflexible,
  keepspaces=true,
  xleftmargin=0.6em,
  xrightmargin=0.6em,
}
\usepackage{fancyhdr}

\fancypagestyle{firstpage}{
  \fancyhf{} 
  \lhead{\includegraphics[height=1.5cm]{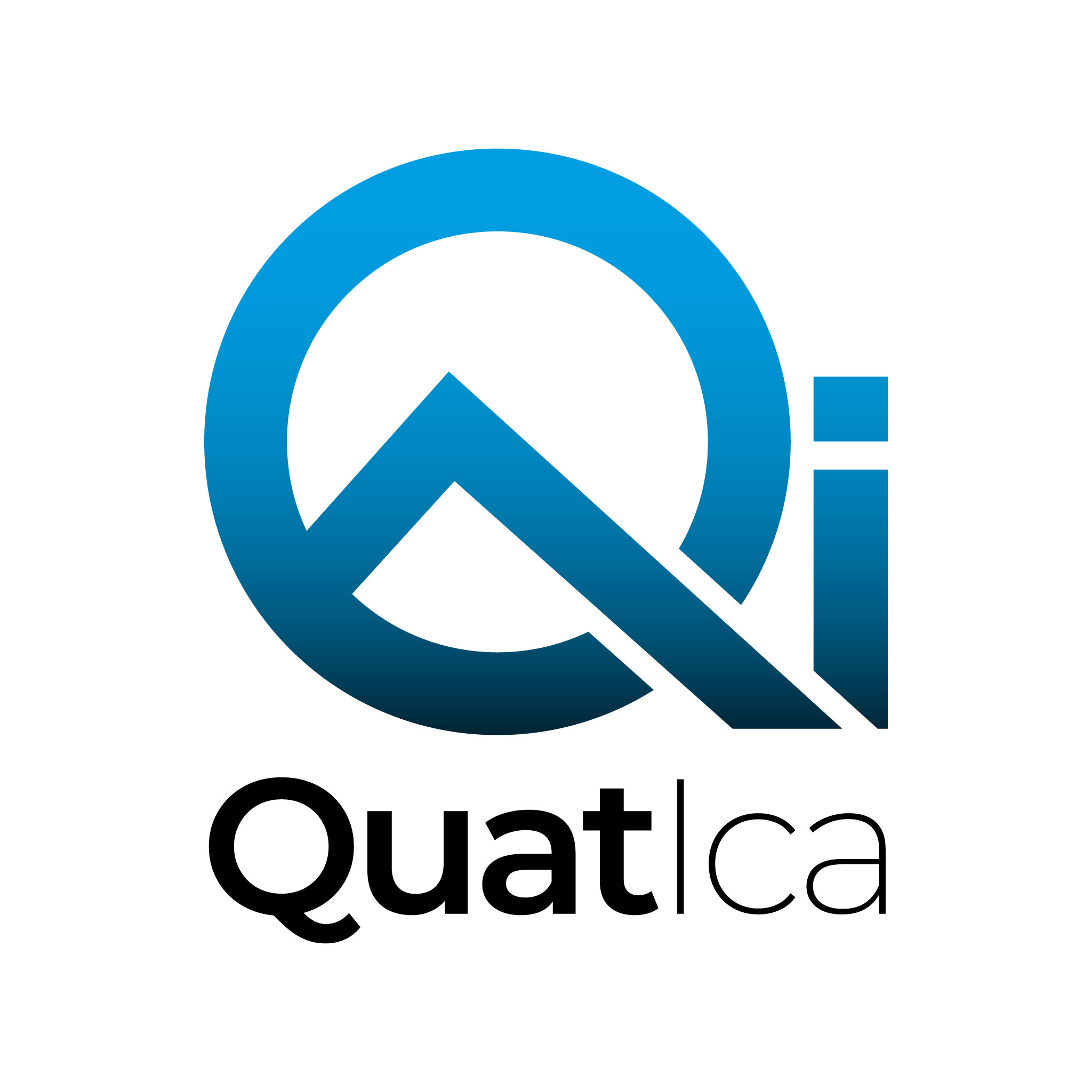}}
  \rhead{\texttt{Quaternion NumerIca (QuatIca)}}
  \cfoot{\thepage}
}

\title{\texttt{QuatIca}: Advanced Numerical Linear Algebra and Optimization for Quaternionic Matrices in Python}

\author[1]{Valentin Leplat}
\author[2]{Salman Ahmadi-Asl}
\author[3]{Junjun Pan}
\author[4]{Henni Ouerdane}
\author[5]{Michael Ng}

\affil[1]{Faculty of Engineering, Innopolis University, Russia\\
\texttt{v.leplat@innopolis.ru}}
\affil[2]{Lab of Machine learning and knowledge representation, Innopolis University, Russia\\
\texttt{s.ahmadiasl@innopolis.ru}}
\affil[3]{Department of Mathematics, Faculty of Science, Hong Kong Baptist University, Hong Kong, China.\\
\texttt{junjpan@hkbu.edu.hk}}
\affil[4]{Center for Engineering Systems and Sciences, Skolkovo Institute of Science and Technology, Moscow, Russia\\
\texttt{h.ouerdane@skoltech.ru}}
\affil[5]{Department of Mathematics, Faculty of Science, Hong Kong Baptist University, Hong Kong, China\\
\texttt{michael-ng@hkbu.edu.hk}}

\date{}

\begin{document}
\maketitle
\thispagestyle{firstpage} 

\begin{abstract}
Quaternion-valued representations provide a convenient way to model coupled multi-channel signals (e.g., RGB imagery, polarization data, vector fields, and multi-detector time series). Yet practical and numerically reliable software support remains far less mature than 
those based on the real/complex setting. Here, we present \texttt{QuatIca}, an open-source Python library for quaternion numerical linear algebra and optimization, designed for both research prototyping and reproducible experimentation. \texttt{QuatIca} provides core quaternion matrix operations and norms; dense decompositions and reductions (QR, LU, Q-SVD, eigendecomposition, Hessenberg/tridiagonal reduction, Cholesky decomposition, and Schur helpers); iterative solvers including quaternion GMRES (with preconditioning) and Newton-Schulz pseudoinverse schemes; and domain-focused routines for signal and image processing such as quaternion Tikhonov restoration. The library also includes \texttt{OptiQ}, which solves quaternion Hermitian semidefinite programs using log-det barrier Newton methods with $\mu$-continuation. We highlight design choices that preserve quaternion structure, and we provide end-to-end demonstrations
including quaternion image deblurring, Lorenz-attractor filtering, and quaternion image completion. 
\texttt{QuatIca} is distributed via PyPI and accompanied by open-source development on GitHub and continuously deployed documentation with runnable tutorials.
\end{abstract}


\section{Introduction}
Quaternion-valued representations offer advantages over traditional methods in the field of multidimensional signal processing. When handling
multi-channel signals whose components are strongly coupled, 
they conserve meaningful cross-channel structure.
Typical examples include RGB imagery encoded as pure quaternions, 
polarized signals, 3D vector fields, and multi-detector time series. In such settings, quaternion algebra provides a compact formalism for preserving inter-channel coupling and phase information. It also introduces distinctive numerical challenges, including noncommutativity, quaternion Hermitian structure, and right/left linearity conventions.
Despite the long history of quaternion linear algebra and its many potential advantages, high-quality, research-oriented numerical software remains underdeveloped compared to the real and complex settings.

\paragraph{Contributions.}
In this article, we introduce \texttt{QuatIca}, a Python library that aims to fully integrate 
quaternion numerical linear algebra into modern scientific computing workflows. The library is designed to be:
\begin{itemize}[leftmargin=1.2em]
\item \emph{Practical}: with simple APIs, copy-paste examples, and end-to-end demonstrations;
\item \emph{Numerically sound}: with explicit attention to quaternion-specific algebraic structure and numerical stability;
\item \emph{Extensible}: with a modular design that supports new solvers and application-specific pipelines, including conic optimization through \texttt{OptiQ}.
\end{itemize}

\paragraph{Availability, dependencies, and reproducibility.}
\texttt{QuatIca} is released as an installable Python package on PyPI\footnote{\url{https://pypi.org/project/quatica/}}
and can be installed via \texttt{pip} (e.g., \texttt{pip install quatica}). Recognizing the growing importance of Python in scientific computing and the need for robust quaternion tools in the Python ecosystem, \texttt{QuatIca} is designed to build on the strong foundations of standard scientific Python. In particular, it relies on widely used and well-tested libraries such as NumPy, SciPy, and Matplotlib, and leverages the \texttt{numpy-quaternion} package for fundamental quaternion arithmetic and data structures. This ensures that higher-level quaternion matrix and tensor routines in \texttt{QuatIca} integrate naturally with existing Python workflows. The full source code is available on GitHub\footnote{\url{https://github.com/vleplat/QuatIca}} together with unit tests, reproducible scripts, and examples. A continuously deployed documentation website (MkDocs)\footnote{\url{https://vleplat.github.io/QuatIca/}} provides API references and runnable tutorials, including links to online Google Colab notebooks for interactive experimentation without local installation.

\paragraph{Positioning with respect to existing software.}
Quaternion numerical computing has been supported mainly through specialized environments, most notably the Quaternion Toolbox for MATLAB (QTFM)\footnote{\url{http://qtfm.sourceforge.net/} (Stephen J. Sangwine and Nicolas Le Bihan).}, while the Python ecosystem has so far offered primarily lower-level quaternion packages. In particular, \texttt{numpy-quaternion} provides a quaternion data type and basic quaternion arithmetic for NumPy arrays, together with several geometric utilities, but it does not aim to provide a full framework for quaternion matrix numerical linear algebra and optimization. To the best of our knowledge, \texttt{QuatIca} is the first open-source scientific Python framework specifically dedicated to quaternion data analysis in this broader sense, combining core quaternion matrix operations with advanced decompositions, iterative solvers, optimization routines, and reproducible application pipelines within a unified research-oriented library.

\paragraph{Article outline.}
In Section~\ref{sec:math}, we provide a summary of quaternion basic definitions and properties 
used throughout the work.
Section~\ref{sec:software} is devoted to the software architecture and main modules. Key algorithms implemented in \texttt{QuatIca} are described in Section~\ref{sec:algorithms}.
In Section~\ref{sec:applications} several practical examples of the application of \texttt{QuatIca} for various problems are given. 
In Section~\ref{sec:repro}, we discuss the testing framework, numerical validation strategy, and benchmark experiments.
The article concludes on current limitations and an outline of future directions in Section~\ref{sec:outlook}.

\section{Quaternion preliminaries}
\label{sec:math}
A quaternion $q\in\HH$ is of the form $q=q_0+q_1\mathtt{i}+q_2\mathtt{j}+q_3\mathtt{k}$, where $q_i$ ($i=0,1,2,3$) are real scalars, and $\mathtt{i},\mathtt{j},\mathtt{k}$ are basis elements satisfying the quaternionic relations: $\mathtt{i}^2=\mathtt{j}^2=\mathtt{k}^2=\mathtt{ijk}=-1$. The conjugate of a quaternion $q$ reads:
$\bar q=q_0-q_1\mathtt{i}-q_2\mathtt{j}-q_3\mathtt{k}$, and its modulus is $|q|=\sqrt{q\bar q}$. For a quaternion matrix $A\in\HH^{m\times n}$, its conjugate-transpose is
$A^H := (\overline{A})^T$. A quaternion matrix is \emph{Hermitian} if $A^H=A$. A semidefinite positive matrix and a positive definite matrix are denoted by $X \succeq 0$ and $X \succ 0$, respectively.
\paragraph{Real inner product on $\Herm_n(\HH)$.}
Let $\Herm_n(\HH)$ denote the set of $n\times n$ quaternion Hermitian matrices.
This set forms a real vector space under addition and multiplication by real scalars.
We use the real inner product
\[
\ip{U}{V} := \Rea\,\tr(U^H V),\qquad U,V\in\Herm_n(\HH),
\]
which is consistent with many Hermitian-structure-preserving algorithms.
The associated Frobenius norm is
\[
\|U\|_F := \sqrt{\ip{U}{U}}.
\]

\paragraph{Adjoint vs Hermitian projection (implementation-critical).}
We distinguish:
\[
A^H\ \text{(adjoint)}\qquad \text{vs}\qquad \herm(A):=\tfrac12(A+A^H)\in\Herm_n(\HH)\ \text{(projection)}.
\]
The latter is used to remove numerical drift and enforce Hermitian constraints.

\paragraph{Unit quaternions and 3D rotations (\texttt{qtraj}).}
In addition to matrix computations, \texttt{QuatIca} includes a module \texttt{qtraj} for quaternion
\emph{rotation trajectories}. The relevant state space is the set of \emph{unit quaternions}
\[
S^3 := \{q\in\HH:\ |q|=1\}\subset\R^4,
\]
which forms a group under quaternion multiplication (with inverse $q^{-1}=\bar q$ for $q\in S^3$). 
A vector $x=(x_1,x_2,x_3)\in\R^3$ is identified with the pure quaternion
$x=x_1\mathtt{i}+x_2\mathtt{j}+x_3\mathtt{k}$ (zero scalar part). 

\paragraph{Rotation matrices and the map $S^3\to SO(3)$.}
The group of proper 3D rotations is
\[
SO(3)=\{R\in\R^{3\times 3}: R^\top R=I,\ \det(R)=1\}.
\]
Given $q\in S^3$ and a pure quaternion $x\in\R^3$, we define
\[
x' = q\,x\,q^{-1}.
\]
One can show that $x'$ is again pure and that the induced linear map $x\mapsto x'$ is a rotation, which thus corresponds to a matrix $R(q)\in SO(3)$ \cite{altmann2005rotations}. 
A key subtlety is the \emph{double cover}: $R(q)=R(-q)$, so each rotation has two representing unit quaternions $\{q,-q\}$. 

\paragraph{Sign consistency (unwrapping) and the dot product in $\R^4$.}
Because of the sign ambiguity, keyframes are typically ``unwrapped'' so that consecutive quaternions stay close on $S^3$. \texttt{qtraj} follows the standard convention: if $\langle q_{i-1},q_i\rangle<0$, then replace
$q_i\leftarrow -q_i$. 
Here $\langle\cdot,\cdot\rangle$ is the Euclidean dot product in $\R^4$:
for $q=(q_0,q_1,q_2,q_3)$ and $r=(r_0,r_1,r_2,r_3)$,
\[
\langle q,r\rangle := q_0 r_0 + q_1 r_1 + q_2 r_2 + q_3 r_3.
\]

\paragraph{Axis--angle representation and quaternion exponential/logarithm.}
Any 3D rotation can be described by an axis \(u\in\R^3\) with \(\|u\|=1\) and a rotation angle \(\theta\).
When \(u\) is identified with the pure quaternion
\[
u=u_1\mathtt{i}+u_2\mathtt{j}+u_3\mathtt{k},
\]
the corresponding unit quaternion is
\[
q=\cos\!\Big(\frac{\theta}{2}\Big)+u\,\sin\!\Big(\frac{\theta}{2}\Big).
\]
Thus, the direction of the imaginary part of \(q\) gives the rotation axis, while its magnitude encodes half of the rotation angle.

The quaternion exponential provides a compact way to build such a unit quaternion from an axis--angle quantity.
Indeed, for a pure quaternion \(v\neq 0\), written as
\[
v=\alpha\,\hat v,
\qquad
\alpha=\|v\|,
\qquad
\|\hat v\|=1,
\]
one has
\[
\exp(v)=\cos(\alpha)+\hat v\,\sin(\alpha).
\]
In particular, if \(v=u\,\theta/2\), then
\[
\exp(v)=\cos\!\Big(\frac{\theta}{2}\Big)+u\,\sin\!\Big(\frac{\theta}{2}\Big),
\]
which is exactly the quaternion representing the rotation of angle \(\theta\) around axis \(u\).

Conversely, the quaternion logarithm recovers this axis--angle information locally.
If a unit quaternion is written as
\[
q=q_0+q_v,
\]
with scalar part \(q_0\) and pure part \(q_v\), then, away from the degenerate case \(q_v=0\),
\[
\log(q)=\alpha\,\frac{q_v}{\|q_v\|},
\qquad
\alpha=\arccos(q_0).
\]
Hence, \(\log(q)\) is a pure quaternion whose direction gives the rotation axis and whose norm equals half of the rotation angle. More broadly, the exponential and logarithm maps provide a convenient local bridge between rotations, represented by unit quaternions on \(S^3\), and ordinary vectors in \(\R^3\). In practice, this makes it possible to build rotations from axis--angle data, compare nearby orientations, and construct smooth interpolation schemes using standard Euclidean tools before mapping the result back to the unit-quaternion manifold. This viewpoint underlies the log--exp interpolation routines implemented in \texttt{qtraj}.\\


\section{Software architecture and design goals}
\label{sec:software}

\subsection{Design principles}
\texttt{QuatIca} is guided by three main design principles:
(i) \emph{clean APIs} for quaternion matrices that integrate naturally with standard scientific Python workflows and libraries (NumPy, SciPy, and Matplotlib; see Table~\ref{tab:quatica_glossary});
(ii) \emph{robust numerics} through explicit structure-preserving operations, including careful handling of quaternion Hermitian structure, inner products, and stable factorization and solver building blocks;
and
(iii) \emph{reproducibility} through runnable demonstrations, regression tests, and standardized diagnostics that support reliable research use and long-term maintenance.

\subsection{Module overview}
\label{subsec:module_overview}
\begin{figure}[t]
\centering
\begin{tikzpicture}[
  font=\footnotesize,
  box/.style={draw=black!55, rounded corners=2mm, inner sep=4pt, align=center, line width=0.4pt},
  big/.style={box, minimum width=0.92\linewidth},
  mid/.style={box, minimum width=0.44\linewidth},
  arr/.style={-{Latex[length=3mm]}, thick, draw=black!70, shorten <=2pt, shorten >=2pt}
]

\node[big, fill=gray!6] (deps) {
  \textbf{Scientific Python foundation}\\
  \texttt{numpy-quaternion} \quad+\quad NumPy \quad+\quad SciPy \quad+\quad Matplotlib
};

\node[big, above=7mm of deps, fill=blue!6] (core) {
  \textbf{Core quaternion matrix layer}\\
  \texttt{utils}: norms and inner products, Hermitian tools,\\
  embeddings (\texttt{Realp}/\texttt{Realr}), matrix operations,\\
  null-space utilities, validation helpers
};

\node[mid, above=8mm of core, xshift=-0.30\linewidth, fill=green!6] (decompL) {
  \textbf{Matrix decompositions}\\
  Q-SVD (classical, embedding)\\
  randomized / pass-efficient Q-SVD\\
  QR (via real embedding)\\
  MaxVol / CUR sampling (experimental)
};

\node[mid, above=8mm of core, xshift=+0.30\linewidth, fill=teal!6] (decompR) {
  \textbf{Reductions \& spectral analysis}\\
  LU (native, pivoted)\\
  Hessenberg reduction\\
  Hermitian tridiagonalization\\
  eigendecomposition / Cholesky / Schur
};

\node[mid, above=8mm of decompL, fill=orange!7] (solvers) {
  \textbf{Solvers and restoration}\\
  \texttt{solver.py}: Q-GMRES, pseudoinverse (NS),\\
  CG, randomized sketch-and-project\\
  \texttt{qslst.py}: quaternion Tikhonov restoration\\
  (FFT / matrix variants)
};

\node[mid, above=8mm of decompR, fill=purple!6] (optiq) {
  \textbf{OptiQ (\texttt{optiQ.py})}\\
  Quaternion Hermitian SDP solver\\
  log-det barrier IPM (\textbf{primal})\\
  ($\mu$-continuation + hat-space whitening)
};

\node[mid, above=10mm of $(solvers.north)!0.5!(optiq.north)$, fill=red!6] (qtraj) {
  \textbf{qtraj (\texttt{qtraj.py})}\\
  Rotation trajectories on $S^3$\\
  SLERP, SQUAD, log--exp
};

\node[big, below=7mm of deps, fill=yellow!10] (io) {
  \textbf{Diagnostics \& plotting}\\
  \texttt{visualization}: standardized figures and saved logs
};

\draw[arr] (deps.north) -- (core.south);

\draw[arr] (core.north) -- (decompL.south);
\draw[arr] (core.north) -- (decompR.south);

\draw[arr] (decompL.north) -- (solvers.south);
\draw[arr] (decompR.north) -- (optiq.south);

\coordinate (split) at ($(decompL.north)!0.5!(decompR.north)+(0,6mm)$);
\draw[arr] (core.north) -- ++(0,3mm) |- (split);
\draw[arr] (split) -| (solvers.south);
\draw[arr] (split) -| (optiq.south);

\draw[arr] (core.north) -- ++(0,3mm) |- (qtraj.south);

\draw[arr] (deps.south) -- (io.north);

\coordinate (routeL) at ($(deps.west)+(-8mm,0)$);
\draw[arr] (core.south) -- ++(0,-3mm) -| (routeL) |- ($(io.north)+(0,2mm)$);

\end{tikzpicture}
\caption{High-level architecture of \texttt{QuatIca}: the core quaternion matrix layer (\texttt{utils}) supports decompositions and spectral reductions, iterative solvers and restoration modules, conic optimization (\texttt{OptiQ}), and rotation trajectories (\texttt{qtraj}), together with standardized diagnostics and plotting.}
\label{fig:architecture}
\end{figure}
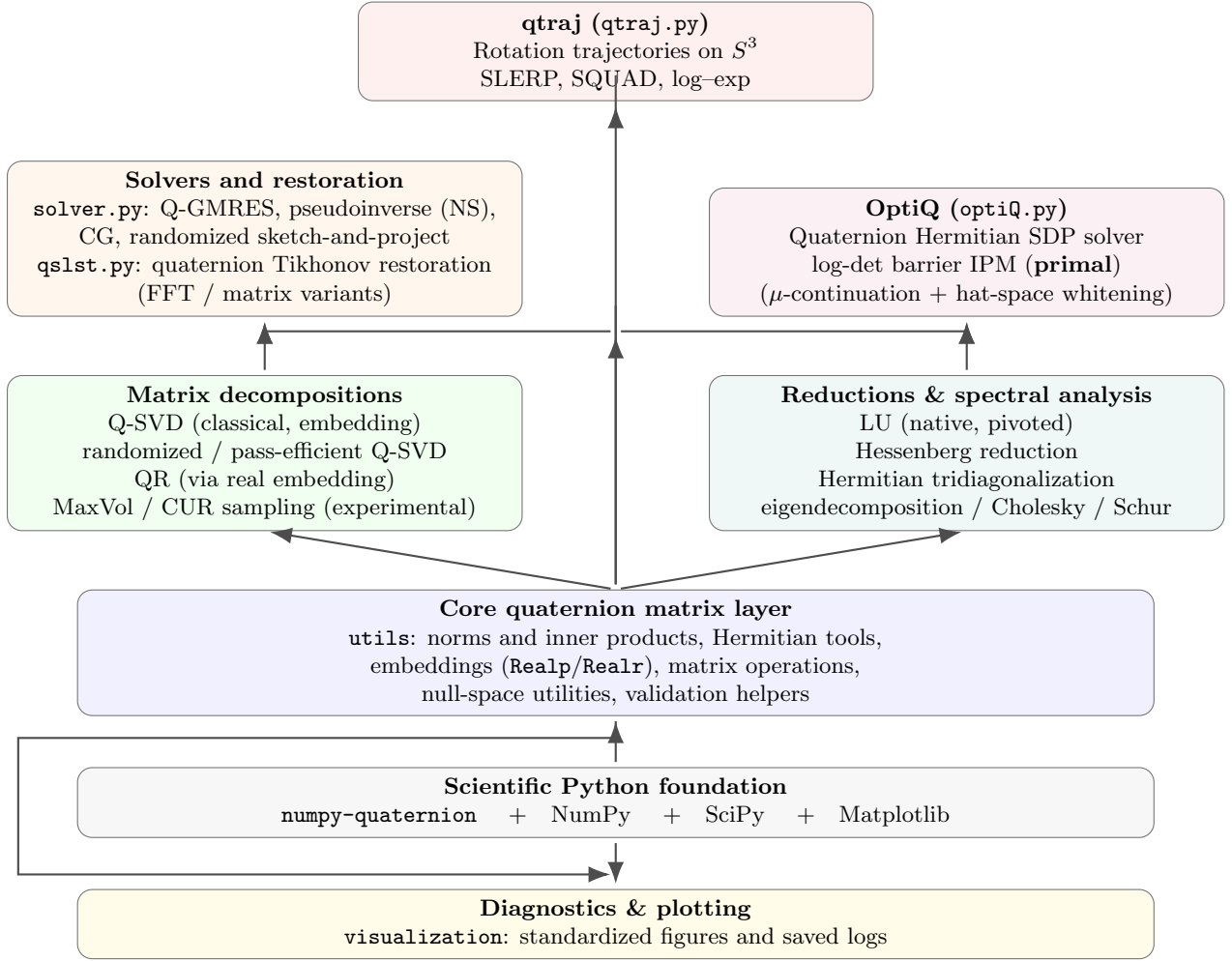

\texttt{QuatIca} follows a layered design. It builds on the standard scientific Python stack for performance and interoperability, while exposing higher-level quaternion routines through a small number of coherent modules. Figure~\ref{fig:architecture} summarizes this organization and highlights how low-level quaternion primitives support decompositions, iterative solvers and restoration methods, conic optimization, and trajectory processing.

For convenience, Table~\ref{tab:quatica_glossary} provides a brief glossary of the main libraries and tools used throughout the \texttt{QuatIca} ecosystem.

\begin{table}[t]
\centering
\caption{Glossary of acronyms and terms for key libraries and tools used by \texttt{QuatIca}.}
\label{tab:quatica_glossary}
\begin{tabularx}{\linewidth}{@{}lX@{}}
\toprule
\textbf{Acronym} & \textbf{Description} \\
\midrule
NumPy & Numerical Python: array programming and core numerical routines. \\
SciPy & Scientific Python: scientific computing routines, including linear algebra and optimization. \\
Matplotlib & Plotting library used for figures and diagnostics. \\
\texttt{numpy-quaternion} & Quaternion dtype and arithmetic for NumPy arrays; the foundation for quaternion data structures in \texttt{QuatIca}. \\
BLAS/LAPACK & Low-level dense linear algebra kernels used indirectly through NumPy/SciPy backends. \\
PyPI & Python Package Index, used to distribute \texttt{QuatIca} as an installable package. \\
\texttt{pip} & Standard Python package installer (e.g., \texttt{pip install quatica}). \\
MkDocs & Static documentation generator used for the online documentation site. \\
Colab & Google Colaboratory: hosted Jupyter notebooks for interactive demos without local setup. \\
GitHub & Source hosting and issue tracking for development and reproducibility artifacts. \\
pytest & Python testing framework used for unit and regression tests. \\
\bottomrule
\end{tabularx}
\end{table}

\paragraph{Core modules.}
At the code level, \texttt{QuatIca} is organized into a small number of modules with clear responsibilities:
\begin{itemize}[leftmargin=1.2em]
\item \texttt{utils.py}: quaternion matrix operations, norms, conversions and embeddings, null-space and kernel utilities, and Hermitian tools.
\item \texttt{decomp/}: Q-SVD (including randomized variants), eigendecomposition, LU and QR, Hessenberg and tridiagonal reductions, Cholesky factorization, and Schur decomposition.
\item \texttt{solver.py}: Newton--Schulz pseudoinverse schemes (order~2 and higher-order); Q-GMRES for quaternion linear systems (optionally LU-preconditioned); and additional pseudoinverse solvers from~\cite{leplat2025iterativemethodscomputingmoorepenrose}, including RSP-Q, hybrid RSP+NS, and CGNE--Q.
\item \texttt{qslst.py}: Quaternion Special Least Squares with Tikhonov regularization, including FFT-specialized and matrix-based restoration variants.
\item \texttt{tensor.py}: experimental quaternion tensor utilities (order~3), unfolding/folding operations, and norms.
\item \texttt{visualization.py}: plotting utilities and standardized diagnostic figures.
\item \texttt{optiQ.py}: \texttt{OptiQ} for quaternion Hermitian semidefinite programs via log-det barrier interior-point methods.
\item \texttt{qtraj.py}: utilities for smooth orientation trajectories on $S^3$ (unit quaternions), including SLERP and SQUAD~\cite{Shoemaker1985,DamKochLillholm1998}, together with a log--exp spline construction and basic kinematic diagnostics~\cite{DamKochLillholm1998,ParkerIbarraOber2023}.
\end{itemize}

\subsection{Reproducible examples and scripts}

Beyond the core API, the repository includes runnable examples, command-line scripts, and notebooks that illustrate the main workflows of the framework.
Several Google Colab notebooks are also linked from the \texttt{PyPI} project page for interactive demonstrations.
These examples serve two complementary purposes: they provide quick-start entry points for new users, and they also act as lightweight end-to-end checks for key applications.

All unit tests can be run from the repository root with
\begin{lstlisting}[language=bash]
pytest -q
\end{lstlisting}
The application scripts can be reproduced either by following the commands documented in the GitHub repository, the framework documentation, and the linked Colab notebooks, or by using the listings provided throughout this article.

\section{Core algorithms implemented in \texttt{QuatIca}}
\label{sec:algorithms}

\subsection{Quaternion decompositions and matrix analysis}

\paragraph{Q-SVD and related decompositions.}
A central goal of \texttt{QuatIca} is to make quaternion matrix analysis as accessible and reliable as its real and complex counterparts, while remaining explicit about quaternion-specific structure. To this end, the library implements the quaternion singular value decomposition (Q-SVD) as a foundational primitive for tasks such as quaternion PCA, low-rank approximation, and denoising, and complements it with standard matrix reductions including Hessenberg form for general matrices and Hermitian tridiagonalization for Hermitian matrices. These decompositions provide numerically stable entry points for eigenvalue computations, Schur-type routines, and Krylov methods used throughout the toolbox. The library is quaternion-native at the interface and algorithmic level, while a small number of routines currently use real-block embeddings as an internal implementation strategy to exploit optimized linear algebra kernels; outputs are subsequently contracted back to quaternion form and checked through systematic validation procedures.


\paragraph{QR decomposition.}
At present, \texttt{QuatIca} computes the QR factorization of a quaternion matrix through the standard $4\times 4$ real-block embedding. For $X\in\HH^{m\times n}$, the routine \texttt{qr\_qua} (currently located in \texttt{quatica/decomp/qsvd.py}) constructs the embedded real matrix $\mathrm{Realp}(X)\in\mathbb{R}^{4m\times 4n}$ using \texttt{real\_expand(X)}, computes a dense QR factorization with SciPy/LAPACK, and contracts the resulting factors back to quaternion form via \texttt{real\_contract}. For tall matrices ($m\ge n$), the routine returns a thin factorization with $Q\in\HH^{m\times n}$ and $R\in\HH^{n\times n}$, where $Q^H Q \approx I_n$ and $X \approx QR$ up to numerical precision. For wide matrices ($m<n$), it returns $Q\in\HH^{m\times m}$ unitary and $R\in\HH^{m\times n}$ upper trapezoidal. This routine serves as a basic orthonormalization primitive and is used in particular within the randomized and pass-efficient Q-SVD pipelines, where stable column orthogonalization is essential. A minimal usage example is shown in Listing~\ref{lst:QR-decomposition}.

\begin{lstlisting}[caption={QR decomposition},label={lst:QR-decomposition}]
from quatica.decomp.qsvd import qr_qua

# QR decomposition of a quaternion matrix
Q, R = qr_qua(X_quat)
# X_quat = Q @ R, where Q has orthonormal columns and R is upper triangular

\end{lstlisting}

\paragraph{LU decomposition.}
\texttt{QuatIca} implements quaternion LU factorization with partial pivoting through the routine \texttt{quaternion\_lu} in \texttt{quatica/decomp/LU.py}, following the logic of the MATLAB QTFM toolbox\footnote{\url{http://qtfm.sourceforge.net/} (Stephen J. Sangwine and Nicolas Le Bihan).}. At each elimination step, the pivot row is selected by maximizing the quaternion modulus in the active column, which improves numerical robustness in the noncommutative setting. The implementation supports rectangular matrices: for $A\in\HH^{m\times n}$, it returns factors $L\in\HH^{m\times N}$ and $U\in\HH^{N\times n}$ with $N=\min(m,n)$. If \texttt{return\_p=True}, the routine returns $(L,U,P)$ such that
\[
PA = LU,
\]
where $P\in\HH^{m\times m}$ is a permutation matrix. If \texttt{return\_p=False}, it returns $(L_{\mathrm{perm}},U)$ such that
\[
A = L_{\mathrm{perm}}U,
\]
where $L_{\mathrm{perm}} = P^T L$ absorbs the row permutations into the lower factor. In \texttt{QuatIca}, this decomposition is also used as a practical preconditioner for iterative solvers such as Q-GMRES, and as a direct-solve baseline in numerical benchmarks. Listing~\ref{lst:LU-decomposition} illustrates a minimal call to \texttt{quaternion\_lu}.

\newpage 
\begin{lstlisting}[caption={LU decomposition},label={lst:LU-decomposition}]
from quatica.decomp import quaternion_lu

# LU decomposition with permutation: P @ A = L @ U
L, U, P = quaternion_lu(A_quat, return_p=True)
# A_quat = (P.T @ L) @ U, where L is lower triangular with unit diagonal
# U is upper triangular, and P is a permutation matrix

# LU decomposition without explicit permutation output
L_perm, U_perm = quaternion_lu(A_quat, return_p=False)
# A_quat = L_perm @ U_perm
# Here L_perm = P.T @ L absorbs the pivoting, so it is generally not lower triangular

\end{lstlisting}


\paragraph{Q-SVD.}
\texttt{QuatIca} provides a classical quaternion singular value decomposition through the routines
\texttt{classical\_qsvd} and \texttt{classical\_qsvd\_full} in \texttt{quatica/decomp/qsvd.py}.
For a matrix $X\in\HH^{m\times n}$, the implementation forms the standard real-block expansion
$\mathrm{Realp}(X)\in\mathbb{R}^{4m\times 4n}$, computes an SVD of the expanded matrix, contracts the
left and right factors back to quaternion form, and recovers the quaternion singular values by retaining
one representative from each fourfold block of repeated singular values induced by the embedding.
The routine \texttt{classical\_qsvd} then truncates the decomposition to a target rank $R$, whereas
\texttt{classical\_qsvd\_full} returns the full factorization. In \texttt{QuatIca}, this Q-SVD serves as a
basic primitive for low-rank approximation, denoising, and spectral diagnostics.

\begin{lstlisting}[caption={Quaternion SVD},label={lst:Q-SVD}]
from quatica.decomp.qsvd import classical_qsvd, classical_qsvd_full

# Truncated Q-SVD for low-rank approximation
U, s, V = classical_qsvd(X_quat, R)
# X_quat \approx U @ diag(s) @ V^H

# Full Q-SVD
U_full, s_full, V_full = classical_qsvd_full(X_quat)
# X_quat = U_full @ \Sigma @ V_full^H

\end{lstlisting}

\paragraph{Randomized Q-SVD and pass-efficient variants.}
For larger matrices, \texttt{QuatIca} includes two approximate low-rank Q-SVD routines:
the randomized method \texttt{rand\_qsvd}, which combines Gaussian sketching with power iterations,
and the pass-efficient method \texttt{pass\_eff\_qsvd}, which reduces data movement by alternating
single multiplications by $X$ or $X^H$ with QR orthonormalization, see~\cite{ahmadiasl2025passefficientrandomizedalgorithmslowrank} for more details.
In the current implementation, the initial sketch is generated as a dense real Gaussian matrix and
embedded into quaternion form through its real component only.
Both routines use \texttt{qr\_qua} as the orthonormalization primitive, then compute a reduced SVD in
real-block form on a small core factor before lifting the resulting singular vectors back to quaternion form.
As in the classical routine, the quaternion singular values are extracted using the same every-fourth rule
associated with the real embedding.
These methods provide practical low-rank approximations that offer a useful accuracy--runtime trade-off
for large-scale matrices.

\newpage 
\begin{lstlisting}[caption={Randomized Q-SVD},label={lst:randomized-Q-SVD}]
from quatica.decomp.qsvd import rand_qsvd

# Randomized low-rank Q-SVD
U, s, V = rand_qsvd(X_quat, R, oversample=10, n_iter=2)
# X_quat \approx U @ diag(s) @ V^H

\end{lstlisting}

\begin{lstlisting}[caption={Pass-efficient randomized Q-SVD},label={lst:pass-efficient-Q-SVD}]
from quatica.decomp.qsvd import pass_eff_qsvd

# Pass-efficient low-rank Q-SVD
U, s, V = pass_eff_qsvd(X_quat, R, oversample=10, n_passes=4)
# X_quat \approx U @ diag(s) @ V^H

\end{lstlisting}

\paragraph{Hessenberg reduction.}
\texttt{QuatIca} provides Hessenberg reduction for general square quaternion matrices through the routine
\texttt{hessenbergize} in \texttt{quatica/decomp/hessenberg.py}.
Given $A\in\HH^{n\times n}$, the method applies a sequence of unitary Householder similarity transformations
to eliminate entries below the first subdiagonal, producing a factorization
\[
H = P A P^H,
\]
where $H$ is upper Hessenberg and $P$ is unitary.
The implementation reuses the quaternion Householder machinery developed for Hermitian tridiagonalization,
and also includes helper routines such as \texttt{is\_hessenberg} and \texttt{check\_hessenberg} to verify
and clean the Hessenberg structure up to numerical tolerance.
This reduction is a standard preprocessing step for QR- and Schur-type eigenvalue methods, and provides a
useful entry point for experimental quaternion spectral routines.

\begin{lstlisting}[caption={Hessenberg reduction},label={lst:hessenberg-reduction}]
from quatica.decomp.hessenberg import hessenbergize, is_hessenberg

# Reduce a square quaternion matrix to Hessenberg form
P, H = hessenbergize(X_quat)

# Check the structure of the reduced matrix
flag = is_hessenberg(H)

\end{lstlisting}

\paragraph{Hermitian tridiagonalization.}
For Hermitian quaternion matrices, \texttt{QuatIca} implements Householder tridiagonalization through the
routine \texttt{tridiagonalize} in \texttt{quatica/decomp/tridiagonalize.py}.
Given $A\in\Herm_n(\HH)$, the algorithm recursively applies unitary reflectors to obtain
\[
B = P A P^H,
\]
where $P$ is unitary and $B$ is tridiagonal.
In accordance with the structure of Hermitian quaternion matrices, the routine post-processes the result by
removing small off-tridiagonal numerical artifacts and projecting the retained entries to their scalar parts,
thereby returning a real tridiagonal matrix with the same eigenvalues as $A$.
This structure-preserving reduction serves as a basic building block for Hermitian spectral routines in the
library.

\begin{lstlisting}[caption={Hermitian tridiagonalization},label={lst:tridiagonalization}]
from quatica.decomp import tridiagonalize

# Tridiagonalize a Hermitian quaternion matrix
P, B = tridiagonalize(A_quat)

# B is real tridiagonal and similar to A_quat
# P is the unitary transformation matrix

\end{lstlisting}

\paragraph{Cholesky decomposition.}
\texttt{QuatIca} provides two complementary Cholesky-based routines for quaternion Hermitian positive definite (HPD) matrices. For dense inputs, \texttt{chol\_quat\_dense} (in \texttt{quatica/decomp/chol.py}) implements a native quaternion Cholesky factorization
\[
A = L L^H,
\]
where \(L\in\HH^{n\times n}\) is lower triangular with real positive diagonal. The routine operates directly in quaternion arithmetic, optionally symmetrizes the input through a \texttt{hermitianize} option, and supports a small diagonal \texttt{jitter} for numerical stabilization. A companion solver, \texttt{solve\_chol\_quat\_dense}, uses forward and backward substitution with the computed factor and supports both single and multiple right-hand sides.

For sparse quaternion matrices, \texttt{QuatIca} provides \texttt{chol\_quat\_sparse}, which relies on a structure-preserving complex embedding and factors the resulting sparse complex matrix with CHOLMOD through \texttt{scikit-sparse}. In this case, the routine returns a factor object exposing linear solves and, when available, a \texttt{logdet()} method, rather than an explicit quaternion triangular factor. This hybrid design allows \texttt{QuatIca} to retain a quaternion-native dense implementation while leveraging mature sparse direct solvers for large structured problems. Additional implementation details are deferred to Appendix~\ref{app:chol_details}.

\begin{lstlisting}[caption={Dense quaternion Cholesky factorization and solve},label={lst:choleskydecomposition}]
from quatica.decomp import chol_quat_dense, solve_chol_quat_dense

# A: quaternion Hermitian positive definite matrix
# b: quaternion right-hand side
L = chol_quat_dense(A, tol=1e-12, hermitianize=True)
x = solve_chol_quat_dense(L, b)
\end{lstlisting}

\paragraph{Null-space computation.}
\texttt{QuatIca} provides numerical null-space utilities through the routines
\texttt{quat\_null\_space}, \texttt{quat\_null\_right}, \texttt{quat\_null\_left}, and the alias
\texttt{quat\_kernel} in \texttt{quatica/utils.py}.
These routines rely on the full classical Q-SVD: given a quaternion matrix \(A\), they call
\texttt{classical\_qsvd\_full}, estimate the numerical rank from the singular values using a relative
threshold, and extract the singular vectors associated with negligible singular values.
More precisely, if \(s_1\) denotes the largest singular value, then \(s_i\) is treated as numerically zero
whenever \(s_i \le \texttt{rtol}\, s_1\).
The right null space is obtained from the trailing columns of \(V\), yielding a basis \(N_r\) such that
\(A N_r \approx 0\), while the left null space is obtained from the trailing columns of \(U\), yielding
\(N_\ell\) such that \(N_\ell^H A \approx 0\).
In the full-rank case, the routines return an empty quaternion matrix with the appropriate number of rows and
zero columns.
These utilities provide a robust way to diagnose rank deficiency and construct null-space constraints in
validation and benchmarking workflows.

\begin{lstlisting}[caption={Quaternion null spaces via Q-SVD},label={lst:null-spaces}]
from quatica.utils import quat_null_right, quat_null_left

# A: quaternion matrix
N_r = quat_null_right(A, rtol=1e-10)  # A @ N_r \approx 0
N_l = quat_null_left(A, rtol=1e-10)   # N_l^H @ A \approx 0
\end{lstlisting}

\subsection{Newton--Schulz pseudoinverse (order-2 and higher order)}
\label{subsec:NS_solver}

For \(A\in\mathbb{H}^{m\times n}\), the Moore--Penrose (MP) pseudoinverse
\(A^\dagger\in\mathbb{H}^{n\times m}\) is the unique matrix satisfying the following four conditions \cite{ben2003generalized,wang2018generalized} 
\begin{align}\label{def:MP-inverse}
AA^\dagger A = A, \qquad
A^\dagger A A^\dagger = A^\dagger, \qquad
(AA^\dagger)^H = AA^\dagger, \qquad
(A^\dagger A)^H = A^\dagger A.
\end{align}

\texttt{QuatIca} implements the quaternion-native Newton--Schulz family for computing inverses and MP pseudoinverses, following the formulation and analysis developed in \cite{leplat2025iterativemethodscomputingmoorepenrose}. These iterations use only quaternion matrix multiplications and adjoints, which makes them attractive for large-scale settings where decomposition-based methods may be too computationally expensive. 

For the square and invertible case, with \(X_0=\alpha A^H\) and \(\alpha\in(0,2/\|A\|_2^2)\), the classical Newton--Schulz iteration is
\[
X_{k+1}=X_k(2I-AX_k),
\]
which yields the quadratic residual relation
\[
\|I-AX_{k+1}\|_2=\|I-AX_k\|_2^2.
\]

For rectangular full-rank matrices, the update must respect the appropriate side convention. In the full column rank case (\(m\ge n\)), one enforces \(X_kA\to I_n\) through the right deviation
\[
F_k = I_n - X_kA
\]
and uses the damped Newton--Schulz step
\[
X_{k+1}=X_k-\gamma\,(X_kA-I_n)\,X_k,\qquad 0<\gamma\le 1.
\]
This gives the exact recurrence
\[
F_{k+1}=(1-\gamma)F_k+\gamma F_k^2,
\]
while \(AX_k\) converges to the range projector \(P=AA^\dagger\). In the full row rank case (\(m\le n\)), one instead enforces \(AX_k\to I_m\) through the left deviation
\[
E_k = I_m - AX_k
\]
and uses
\[
X_{k+1}=X_k-\gamma\,X_k\,(AX_k-I_m),\qquad 0<\gamma\le 1,
\]
so that
\[
E_{k+1}=(1-\gamma)E_k+\gamma E_k^2,
\]
and \(X_kA\) converges to \(Q=A^\dagger A\). 

To accelerate convergence, \texttt{QuatIca} also includes higher-order hyperpower variants. For order \(p\ge 2\), the linear correction is replaced by a truncated Neumann series:
\[
\text{(left)}\qquad X_{k+1}=X_k\sum_{i=0}^{p-1}E_k^{\,i},
\qquad
\text{(right)}\qquad X_{k+1}=\sum_{i=0}^{p-1}F_k^{\,i}X_k.
\]
These updates satisfy the exact residual recurrences
\[
E_{k+1}=E_k^{\,p},\qquad F_{k+1}=F_k^{\,p},
\]
and therefore achieve local order-\(p\) convergence under the same side-consistent initialization framework. Full derivations, convergence guarantees, and cost-saving polynomial factorizations are given in \cite{leplat2025iterativemethodscomputingmoorepenrose}. 

\begin{lstlisting}[caption={Newton--Schulz pseudoinverse in \texttt{QuatIca}},label={lst:pseudoinverse_demo}]
from quatica.solver import NewtonSchulzPseudoinverse, HigherOrderNewtonSchulzPseudoinverse

# A is a quaternion matrix
ns = NewtonSchulzPseudoinverse(gamma=0.5, max_iter=200, tol=1e-8)
A_pinv_ns, res_ns, _ = ns.compute(A)

hon = HigherOrderNewtonSchulzPseudoinverse(max_iter=100, tol=1e-8)
A_pinv_hon, res_hon, _ = hon.compute(A)
\end{lstlisting}

The practical performance of these order-2 and higher-order Newton--Schulz implementations, including accuracy and runtime comparisons, is assessed in Section~\ref{subsec:ns_benchmark}.

\subsection{Quaternion GMRES and preconditioning}\label{subsec:QGMRES}

\texttt{QuatIca} includes a quaternion GMRES solver for linear systems \(Ax=b\), together with optional LU-based left preconditioning. The solver is designed for quaternion-valued inputs and outputs, while internally relying on a component-wise implementation of the core iteration. This provides a practical way to solve quaternion linear systems without requiring the user to explicitly form real block embeddings, while still allowing comparison with embedded alternatives.

The Generalized Minimal Residual Method (GMRES) \cite{Saad1986GMRES} is a standard Krylov subspace method for large, sparse, and generally non-Hermitian linear systems. Its extension to the quaternion setting, often referred to as Quaternion GMRES (QGMRES)~\cite{jia2021structure,csimcsek2024block}, addresses systems of the form
\[
Ax=b,
\]
where \(A\in\mathbb{H}^{n\times n}\) and \(x,b\in\mathbb{H}^n\).

As in the classical case, QGMRES seeks an approximate solution \(x_k\) whose residual is minimized over a growing quaternion Krylov subspace generated from the initial residual \(r_0=b-Ax_0\):
\[
\mathcal{K}_k(A,r_0)=\mathrm{span}\{r_0,Ar_0,A^2r_0,\dots,A^{k-1}r_0\}.
\]
Because quaternion multiplication is noncommutative, the Arnoldi iteration 
and orthogonality relations must be handled with care. Existing quaternion GMRES formulations either rely on structure-preserving real/complex representations or on direct quaternion orthogonalization procedures \cite{jia2021structure,csimcsek2024block}. The current \texttt{QuatIca} implementation follows a component-wise formulation of QGMRES while retaining a quaternion-valued user interface.

To improve convergence, \texttt{QuatIca} currently supports left LU preconditioning. More precisely, if \(M\approx A\) is a quaternion LU factorization-based preconditioner, the solver applies GMRES to the transformed system
\[
M^{-1}Ax=M^{-1}b.
\]
In the present implementation, \(M\) is obtained from a full quaternion LU factorization with partial pivoting, and the action of \(M^{-1}\) is applied through quaternion forward and backward triangular solves. This preconditioned formulation is particularly useful in the benchmark section~\ref{subsec:qgmres_benchmark}, where it is compared with the unpreconditioned solver.

\texttt{QuatIca} provides the class \texttt{QGMRESSolver}, which returns both the approximate solution and a dictionary containing iteration counts, residual information, and convergence diagnostics. A minimal usage example is shown in Listing~\ref{lst:gmres_demo}.

\newpage
\begin{lstlisting}[caption={Quaternion GMRES in \texttt{QuatIca}},label={lst:gmres_demo}]
from quatica.solver import QGMRESSolver

# Unpreconditioned QGMRES
solver = QGMRESSolver(tol=1e-8, max_iter=100)
x, info = solver.solve(A, b)

# Left-LU-preconditioned QGMRES
solver_prec = QGMRESSolver(tol=1e-8, max_iter=100, preconditioner="left_lu")
x_prec, info_prec = solver_prec.solve(A, b)
\end{lstlisting}

\subsection{Quaternion Tikhonov restoration (QSLST)}\label{subsec:QSLST}

To address ill-posed inverse problems such as color image deblurring, \texttt{QuatIca} implements Quaternion Special Least Squares with Tikhonov regularization (QSLST), following \cite{Fei2025}. Let \(X\) denote the unknown quaternion-valued image and \(B\) the observed degraded image, both of spatial size \(H\times W\). At the image level, the restoration problem is
\[
\min_X \|\mathcal{A}(X)-B\|_F^2 + \lambda \|X\|_F^2,
\]
where \(\mathcal{A}\) is a real linear degradation operator acting on images and \(\lambda>0\) is the Tikhonov regularization parameter. This quadratic penalty stabilizes the inversion and balances data fidelity against the energy of the restored solution.

For RGB data, a standard representation is to encode the image as a pure quaternion field,
\[
X(u,v)=X_r(u,v)\,\mathtt{i}+X_g(u,v)\,\mathtt{j}+X_b(u,v)\,\mathtt{k},
\]
although the implementation also supports more general quaternion-valued images represented componentwise. 
If the observation model is
\[
B = h \ast X + \mathcal N,
\]
where \(\ast\) denotes the two-dimensional discrete convolution of the real point-spread function \(h\) with the quaternion-valued image \(X\), applied componentwise, then \(\mathcal{A}\) is the corresponding blur operator.
The associated normal equation at the image level is
\[
(\mathcal{A}^H \mathcal{A} + \lambda I)X = \mathcal{A}^H B.
\]
Equivalently, after vectorization with \(N=HW\), one obtains the real matrix formulation
\[
(A^\top A + \lambda I)x = A^\top b,
\qquad
x=\mathrm{vec}(X),\quad b=\mathrm{vec}(B),\quad A\in\mathbb{R}^{N\times N}.
\]
Because the degradation operator is real, the same matrix acts independently on each quaternion component.

When the blur is induced by convolution with a real PSF under periodic boundary conditions, the operator is block-circulant with circulant blocks and is diagonalized by the two-dimensional FFT. Writing \(\widehat h=\mathrm{FFT2}(h_{\text{centered/padded}})\) and \(\widehat B=\mathrm{FFT2}(B)\), the restored image is obtained frequency-by-frequency as
\[
\widehat X=\frac{\overline{\widehat h}\,\widehat B}{|\widehat h|^2+\lambda},
\qquad
X=\mathrm{IFFT2}(\widehat X).
\]
Since the blur operator is real, the same scalar frequency filter applies independently to each quaternion component. This yields an efficient FFT-based implementation for convolutional restoration.

\texttt{QuatIca} provides two complementary routines in \texttt{qslst.py}. The routine \texttt{qslst\_restore\_fft} implements the FFT-specialized solution for convolutional blur with periodic boundary conditions, while \texttt{qslst\_restore\_matrix} provides a matrix-based reference implementation for a generic real blur matrix \(A\in\mathbb{R}^{N\times N}\) by forming
\[
T=A^\top A+\lambda I,
\]
and applying the pseudoinverse of \(T\) component-wise to \(A^\top b\). The module also includes utilities to construct Gaussian and motion PSFs, apply blur, add white Gaussian noise, and evaluate restoration quality through PSNR and relative error metrics. An end-to-end quaternion image deblurring and denoising demonstration, including a comparison with a Newton--Schulz-based restoration pipeline from \cite{leplat2025iterativemethodscomputingmoorepenrose}, is presented later in Section~\ref{subsec:image_deblurring}.

\begin{lstlisting}[caption={Quaternion Tikhonov restoration in \texttt{QuatIca}},label={lst:tikhonov_demo}]
from quatica.qslst import qslst_restore_fft, qslst_restore_matrix

# Bq: observed quaternion image as (H, W, 4) float array
# psf: real point-spread function (2D)
# A_mat: real blur matrix of shape (H*W, H*W)
# lam: Tikhonov regularization parameter
X_fft = qslst_restore_fft(Bq, psf, lam, boundary="periodic")
X_mat = qslst_restore_matrix(Bq, A_mat, lam)
\end{lstlisting}

\subsection{Quaternion rotation trajectories (\texttt{qtraj})}
\label{subsec:qtraj}

\paragraph{Problem setup.}
Given orientation keyframes $q_0,\dots,q_n\in S^3$ at times $0=t_0<\cdots<t_n=1$, \texttt{qtraj} constructs
an interpolating trajectory $q:[0,1]\to S^3$ such that $q(t_i)=q_i$. Since infinitely many curves can satisfy
these constraints, \texttt{qtraj} focuses on practical families of interpolants that are easy to compute and
that yield smooth transitions between successive orientations (geodesic segment models and spline-like models).

\paragraph{SLERP and piecewise SLERP.}
For two keyframes $q_0,q_1\in S^3$ (after enforcing sign consistency), define $\theta\in[0,\pi]$ by
$\cos\theta=\langle q_0,q_1\rangle$. The spherical linear interpolation (SLERP) is
\[
\mathrm{slerp}(q_0,q_1;t)=
\frac{\sin((1-t)\theta)}{\sin\theta}\,q_0+\frac{\sin(t\theta)}{\sin\theta}\,q_1,\qquad t\in[0,1].
\]
On a single interval, SLERP follows the shortest great-circle arc on $S^3$ and is smooth in $t$.

A common baseline for multi-keyframe data is piecewise SLERP (apply SLERP on each interval $[t_i,t_{i+1}]$).
This guarantees $C^0$ continuity (the curve hits all keyframes) but is generally not $C^1$ at knot points, which can produce visible ``corners'' in the motion due to abrupt changes in angular velocity. This construction is the classical quaternion geodesic interpolant used in computer graphics and animation; see Shoemaker~\cite{Shoemaker1985} and the detailed exposition in Dam \emph{et al.}~\cite{DamKochLillholm1998}.

\paragraph{Kinematics: angular velocity on \(S^3\).}
For a differentiable trajectory \(q(t)\in S^3\), a common angular-velocity convention is
\[
\omega(t)=2\,\dot q(t)\,q(t)^{-1}\in\operatorname{Im}(\HH)\simeq\R^3,
\]
which is a pure quaternion and can be identified with a vector in \(\R^3\).

In discrete time, a standard midpoint approximation from sampled quaternions is
\[
\omega(t_{k+\frac12}) \approx \frac{2}{\Delta t}\,\log\!\big(q(t_{k+1})\,q(t_k)^{-1}\big),
\qquad
t_{k+\frac12}=\frac{t_k+t_{k+1}}{2}.
\]
We use this quantity in the demonstrations to visualize velocity transitions and compare the smoothness of different interpolation schemes.

\paragraph{SQUAD (Shoemaker’s quaternion spline).}
To obtain smoother multi-keyframe trajectories while retaining efficiency, \texttt{qtraj} implements SQUAD, the classical spline-like quaternion construction introduced by Shoemaker~\cite{Shoemaker1985}; see also Dam \emph{et al.}~\cite{DamKochLillholm1998} for a detailed treatment.
It introduces intermediate control quaternions $a_i$ computed from neighbors:
\[
a_i
= q_i\,\exp\!\left(
-\frac14\Big(\log(q_i^{-1}q_{i-1})+\log(q_i^{-1}q_{i+1})\Big)
\right),\qquad i=1,\dots,n-1,
\]
with endpoint choices such as $a_0=q_0$ and $a_n=q_n$.
On each segment $[t_i,t_{i+1}]$, letting $u=(t-t_i)/(t_{i+1}-t_i)$, the SQUAD curve is
\[
\mathrm{squad}(t)=
\mathrm{slerp}\!~\Big(~
\mathrm{slerp}(q_i,q_{i+1};u),\ \mathrm{slerp}(a_i,a_{i+1};u),\ 2u(1-u)
\Big).
\]

In practice, SQUAD preserves the keyframes exactly and typically produces much smoother transitions than
piecewise SLERP.

\paragraph{Log-Exp interpolation (log-Euclidean spline on $S^3$).}
\texttt{qtraj} also provides a practical spline-like construction based on the quaternion logarithm and exponential maps. This follows the standard idea of mapping nearby orientations to a logarithmic chart, interpolating there with Euclidean spline tools, and mapping back to \(S^3\); see, for example, Dam \emph{et al.}~\cite{DamKochLillholm1998} for quaternion interpolation background and Parker, Ibarra, and Ober~\cite{ParkerIbarraOber2023} for a recent explicit logarithm-based interpolation framework.
Choose a reference quaternion $q_{\mathrm{ref}}$ (commonly $q_0$), map each keyframe to a vector in $\R^3$ by
\[
p_i := \log\!\big(q_{\mathrm{ref}}^{-1}q_i\big)\in\R^3,
\]
fit a cubic spline (or other smooth interpolant) $p(t)$ through $(t_i,p_i)$ in $\R^3$, and map back to $S^3$ via
\[
q(t)=q_{\mathrm{ref}}\,\exp\!\big(p(t)\big).
\]
This \emph{log--exp} approach is simple to implement and often yields very smooth motion; in practice, it works best
when consecutive rotations are not close to the branch cut of the logarithm, which is handled robustly by the
sign-consistency step described in Section~\ref{sec:math}.

\paragraph{Demonstration.}
Section~\ref{sec:applications} presents a dedicated demonstration comparing piecewise SLERP, SQUAD, and log--exp
interpolation on representative rotation sequences, and illustrates the resulting trajectory smoothness through
angular-velocity profiles.

\subsection{OptiQ: Quaternion Hermitian SDP via log-det barrier interior-point methods}
\label{sec:optiq}

\subsubsection{Quaternion Hermitian SDP}
OptiQ targets quaternion Hermitian semidefinite programs (SDPs) of the form
\begin{equation}
\label{eq:optiq_primal_sdp}
\min_{X\in\Herm_n(\HH)}\ \ip{C}{X}
\quad\text{s.t.}\quad A(X)=b,\quad X\succeq 0,
\end{equation}
where the linear operator $A:\Herm_n(\HH)\to\R^m$ is defined through Hermitian measurement matrices
$H_1,H_2,\dots,H_m\in\Herm_n(\HH)$ as
\[
(A(X))_i := \ip{H_i}{X}=\Re\,\tr(H_i^H X),
\qquad i=1,2,\dots,m.
\]
All matrices that should belong to $\Herm_n(\HH)$ are explicitly enforced to be Hermitian in implementation
via the projection $\herm(A)=\tfrac12(A+A^H)$.

\subsubsection{Equality-constrained log-det barrier formulation}
Instead of solving \eqref{eq:optiq_primal_sdp} directly with a primal--dual conic IPM,
OptiQ uses a log-det barrier for the cone constraint \(X\succeq 0\).
For a fixed \(\mu>0\), we solve the equality-constrained barrier subproblem
\begin{equation}
\label{eq:optiq_barrier_subproblem}
\min_{X\in\Herm_n(\HH),\ X\succ 0}\ \ip{C}{X} - \mu \log\det(X)
\quad\text{s.t.}\quad A(X)=b.
\end{equation}

For \(X\in\Herm_n(\HH)\) with \(X\succ0\), the eigenvalues of \(X\) are positive real numbers.
Throughout \texttt{OptiQ}, we define
\[
\log\det(X):=\sum_{i=1}^n \log \lambda_i(X),
\]
where \(\lambda_1(X),\dots,\lambda_n(X)>0\) are the eigenvalues of \(X\).
Equivalently, this agrees with the log-determinant induced by the standard real embedding of
quaternion Hermitian positive definite matrices.

With this convention, and using the real inner product on \(\Herm_n(\HH)\) to identify differentials with gradients, the barrier gradient and Hessian take the familiar form
\[
\nabla \big(-\mu\log\det(X)\big)=-\mu X^{-1},
\qquad
D^2\!\big(-\mu\log\det(X)\big)[W]=\mu X^{-1}WX^{-1},
\quad W\in\Herm_n(\HH).
\]
Here \(W \in \Herm_n(\HH)\) denotes an arbitrary Hermitian perturbation direction, and the second formula is the Hessian viewed as a linear operator acting on \(W\).
This is the convention used throughout the Newton and path-following routines in \texttt{OptiQ}.

The barrier problem~\eqref{eq:optiq_barrier_subproblem} is solved by Newton steps applied to its Karush--Kuhn--Tucker (KKT) optimality system. A path-following (continuation) strategy then decreases \(\mu\) geometrically, using the solution at one barrier level to initialize the next.

\paragraph{Two solver modes.}
OptiQ exposes two closely related routines:
\begin{itemize}[leftmargin=1.2em]
\item a \emph{fixed-\(\mu\) Newton solver} (one solve of \eqref{eq:optiq_barrier_subproblem});
\item a \emph{barrier path-following solver} that repeats fixed-\(\mu\) solves while decreasing \(\mu\)
until a target gap floor \(\varepsilon_{\rm gap}\).
\end{itemize}
Implementation details for the Newton system, Schur complement, and line-search are given in
Appendix~\ref{app:optiq_details}.

\subsubsection{Hat-space constraint whitening (Gram--Cholesky)}
To improve conditioning, OptiQ applies a linear change of coordinates in $\R^m$ based on the Gram matrix
\[
G_{ij}=\ip{H_i}{H_j}\in\R,\qquad G\in\R^{m\times m}.
\]
Assuming the constraints are independent, $G\succ 0$ and admits a Cholesky factorization $G=R^T R$.
We then define the \emph{hat operators}
\begin{equation}
\label{eq:optiq_hat_ops}
\hat A(X) := R^{-T}A(X),\quad
\hat A^*(\hat y):=A^*(R^{-1}\hat y),\quad
\hat b := R^{-T}b.
\end{equation}
A key property is
\begin{equation}
\label{eq:optiq_hat_isometry}
\hat A\hat A^* = I_m,
\end{equation}
so that hat-space residuals are scale-normalized. Hat-space also yields a simple equality correction
(projection onto $\hat A(X)=\hat b$):
\begin{equation}
\label{eq:optiq_equality_snap}
X \leftarrow X + \hat A^*(\hat b-\hat A(X)),
\end{equation}
which is exact whenever \eqref{eq:optiq_hat_isometry} holds. This correction is used throughout the
Newton iterations to control feasibility drift.

\subsubsection{Algorithmic summary}

Our log-det barrier path-following method is summarized in Algorithm~\ref{alg:optiq_barrier_path}.

\begin{algorithm}[ht!]
\caption{Log-det barrier path-following in hat space (OptiQ)}
\label{alg:optiq_barrier_path}
\begin{algorithmic}[1]
\STATE Input: Hermitian constraints $H_i$, vector $b$, cost $C$, initial $X_0\succ 0$, $\mu_0$,
decay $\beta_\mu\in(0,1)$, gap floor $\varepsilon_{\rm gap}$.
\STATE Build hat operators $(\hat A,\hat A^*,\hat b)$ via Gram--Cholesky.
\STATE Initialize $(X,\hat y)\leftarrow (X_0,0)$ and $\mu\leftarrow \mu_0$.
\FOR{stages $s=0,1,\dots$}
  \STATE Warm-start multiplier (optional): $\hat y \leftarrow \hat A(\mu X^{-1}-C)$.
  \STATE Approximately solve the barrier KKT system at current $\mu$ by Newton steps:
  \STATE \hspace{1em} form residuals $\hat r_p=\hat b-\hat A(X)$ and $\hat r_d=C+\hat A^*(\hat y)-\mu X^{-1}$;
  \STATE \hspace{1em} compute Newton direction via Schur complement (Appendix~\ref{app:optiq_details});
  \STATE \hspace{1em} apply fraction-to-boundary and Armijo backtracking;
  \STATE \hspace{1em} apply a damped hat-space feasibility correction toward \(\hat A(X)=\hat b\), preserving \(X\succ0\).
  \IF{$\mu\le \varepsilon_{\rm gap}$} \STATE break \ENDIF
  \STATE Update $\mu \leftarrow \max(\beta_\mu\mu,\varepsilon_{\rm gap})$.
\ENDFOR
\STATE Output: $X(\mu)$ and the associated slack $S(\mu)=\mu X(\mu)^{-1}$.
\end{algorithmic}
\end{algorithm}

\subsubsection{Quick-start: solving a quaternion Hermitian SDP with OptiQ}
Listing~\ref{lst:optiq_quickstart} shows a minimal example using the log-det barrier path solver.

\begin{lstlisting}[caption={OptiQ quick-start (barrier path-following).},label={lst:optiq_quickstart}]
import numpy as np
from quatica import solve_barrier_path  # \mu-continuation logdet barrier IPM

# H_list: list of Hermitian quaternion constraint matrices (n x n)
# b:      real constraint vector (m,)
# C:      Hermitian quaternion cost matrix (n x n)
# X0:     strictly feasible warm start (n x n), X0 $\succ$ 0

out = solve_barrier_path(
    H_list, b, C,
    X0=X0,
    mu0=1e-1,
    beta_mu=0.5,
    eps_gap=1e-10,
    tol=1e-10,
    max_newton=50,
)

X = out["X"]      # primal PSD matrix
S = out["S"]      # recovered slack: S = \mu X^{-1}
yhat = out["yhat"]

print("objective =", out["obj"])
print("||rp_hat|| =", np.linalg.norm(out["rp_hat"]))
\end{lstlisting}

\paragraph{Demonstration.}
Section~\ref{sec:repro} reports an end-to-end validation of \texttt{OptiQ} using a certified known-optimum SDP
instance, illustrating the expected barrier-path behavior and numerical accuracy.

\section{Applications}
\label{sec:applications}

\subsection{Quaternion image completion}
\label{subsec:img_completion}

Image completion is a fundamental task in computer vision: visual data with missing pixels can often be reconstructed reliably from the observed ones when the underlying structure is sufficiently low-dimensional.
A common computational primitive in many completion pipelines is the repeated computation of a low-rank approximation of intermediate iterates.
In this section, we rely on cross approximation (CUR), which samples a subset of columns and rows to obtain an accurate low-rank model at significantly reduced cost.

To preserve cross-channel coupling, we represent an RGB image as a quaternion-valued array by encoding each pixel as a pure quaternion (zero real part),
\[
M_{ij} \;=\; r_{ij}\mathtt{i} + g_{ij}\mathtt{j} + b_{ij}\mathtt{k}\in\HH,
\]
where \((r_{ij},g_{ij},b_{ij})\) are the red, green, and blue intensities.
Setting the scalar part to zero is a standard convention in quaternion color imaging: it ensures that the three imaginary components correspond exactly to the three color channels, while keeping the representation compact and algebraically consistent.

The completion scheme in \cite{ahmadi2023fast,leplat2025iterativemethodscomputingmoorepenrose} iterates
\begin{align}
X^{(n)} &\leftarrow \mathcal{L}\!\big(Y^{(n)}\big), \label{Step1}\\
Y^{(n+1)} &\leftarrow \Omega \oast M + \big(\mathbf{1}-\Omega\big)\oast X^{(n)}, \label{Step2}
\end{align}
where \(M\) is the observed incomplete quaternion-valued image, \(\Omega\in\{0,1\}^{m\times n}\) is the sampling mask
(\(\Omega_{ij}=1\) if \((i,j)\) is observed), \(\mathbf{1}\) is the all-ones matrix, and \(\oast\) denotes the Hadamard product.
The operator \(\mathcal{L}\) computes a low-rank cross approximation in CUR form,
\[
\mathcal{L}(Y)=CUR,\qquad U = C^{\dagger} Y R^{\dagger},
\]
where \(C\) and \(R\) are sampled column and row submatrices of \(Y\), and \((\cdot)^{\dagger}\) denotes the MP pseudoinverse over \(\HH\).
In the current implementation, the required quaternion pseudoinverses are computed through a Newton--Schulz iteration as in \cite{leplat2025iterativemethodscomputingmoorepenrose}. 
The demo script uses a simple randomized row/column sampling strategy inside the CUR step\footnote{More advanced row/column selection strategies could also be incorporated, for example leverage-score-based sampling or pivot-based selection rules.} and applies a light Gaussian denoising step between outer iterations.

\paragraph{Demonstration script.}
The \texttt{QuatIca} repository includes a reproducible demonstration script for quaternion image completion.
In its current default configuration, the script loads a color image, resizes it to \(256\times256\), masks \(70\%\) of the pixels, performs 100 completion iterations, and uses rank parameters \(R_1=R_2=60\).
Running the demo produces:
(i) the original image,
(ii) the masked image,
(iii) the reconstructed image,
(iv) the evolution of the PSNR over iterations,
and
(v) summary information including the final PSNR, the rank parameters, the maximum number of internal Newton--Schulz iterations used in the pseudoinverse subroutines, and the total runtime.
To reproduce the results, one can run the code provided in Listing~\ref{lst:image_completion_demo}.

To the best of our knowledge, no official reference implementation of the quaternion CUR completion pipeline of \cite{wu2025efficient} is currently distributed alongside the article, and \texttt{QuatIca} therefore provides a practical open-source implementation of this type of workflow. 

\begin{lstlisting}[language=bash,caption={Running the quaternion image completion demo.},label={lst:image_completion_demo}]
# From the repository root
python run_analysis.py image_completion
\end{lstlisting}

\begin{figure}[ht!]
  \centering
  \includegraphics[width=\linewidth]{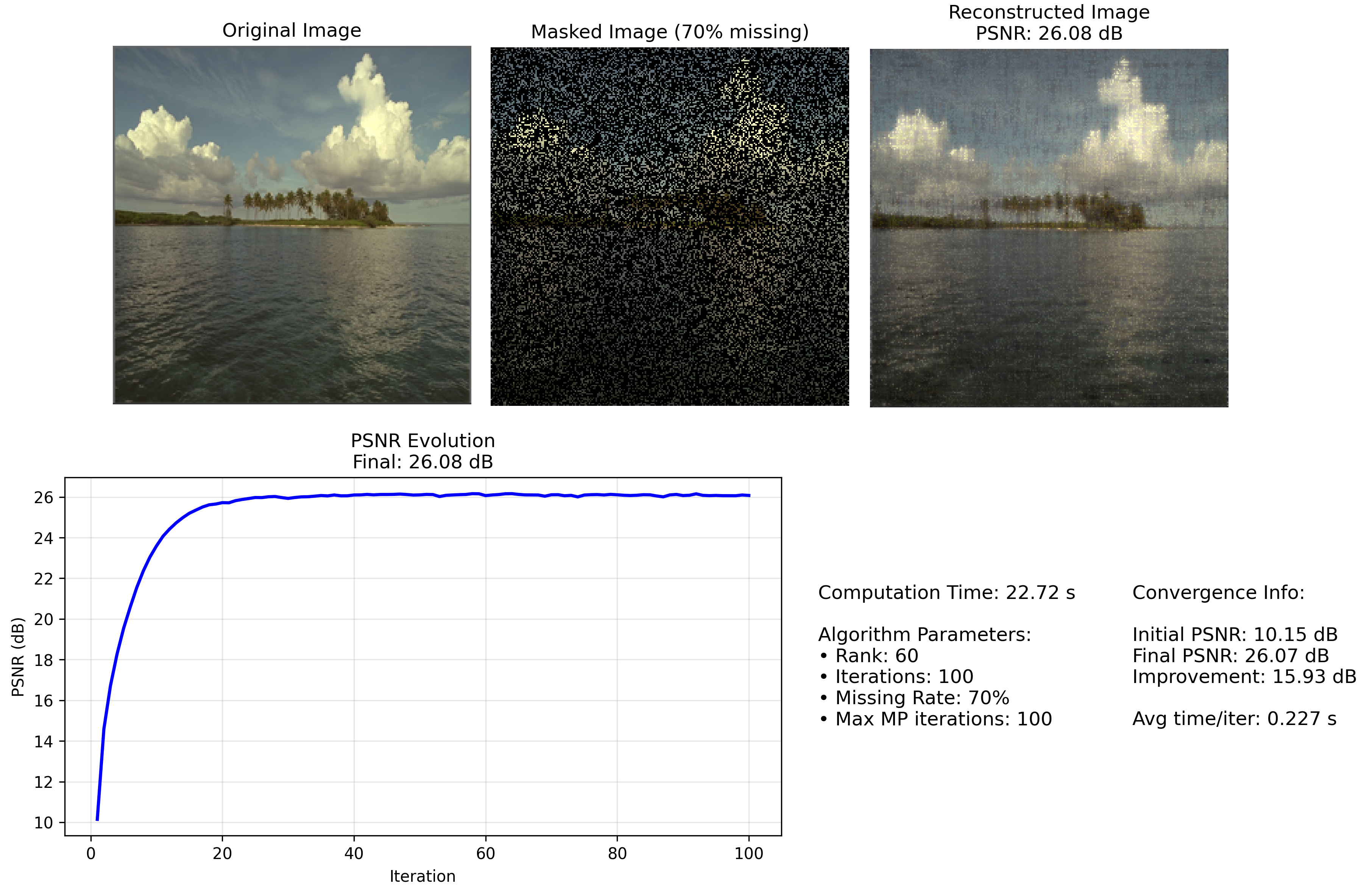}
  \caption{Quaternion image completion demo with \(70\%\) missing pixels. The output summarizes the original image, the masked observation, the final reconstruction, and the evolution of the PSNR across iterations, together with key summary information such as the rank parameters, the pseudoinverse iteration budget, and the total runtime.}
  \label{fig:quat_image_completion_demo}
\end{figure}

\subsection{Quaternion image deblurring}
\label{subsec:image_deblurring}

As in Section~\ref{subsec:img_completion}, RGB images are represented as quaternion-valued arrays, with each pixel encoded as a pure quaternion
\[
p = r\,\mathtt{i} + g\,\mathtt{j} + b\,\mathtt{k}.
\]
This representation preserves the coupling between color channels and provides a natural setting for quaternion-valued restoration methods.

In the deblurring experiments, the observed image is modeled as
\[
B = h \ast X + \mathcal N,
\]
where \(X\) is the unknown clean quaternion image, \(h\) is a real point-spread function (PSF), and \(\mathcal N\) denotes additive noise.
Equivalently, if \(\mathcal A\) denotes the blur operator associated with \(h\), the model can be written as
\[
B = \mathcal A(X) + \mathcal N.
\]
We consider periodic boundary conditions, so that the blur operator is diagonalized by the two-dimensional FFT.
This makes it possible to compare two restoration pipelines under the same degradation model:
(i) the FFT-specialized quaternion Tikhonov solver QSLST implemented in \texttt{qslst.py}, following \cite{Fei2025} and discussed in Section~\ref{subsec:QSLST}, and
(ii) an FFT-based Newton--Schulz inverse iteration for the regularized frequency-domain filter, inspired by the pseudoinverse framework of \cite{leplat2025iterativemethodscomputingmoorepenrose}.

The goal of this demo is not only to produce visually plausible restorations, but also to compare the methods quantitatively through PSNR, SSIM, relative reconstruction error, and runtime.
In the benchmark setting used in this article, both restoration pipelines are run under the same experimental conditions: the same Gaussian blur kernel, the same noise level, the same FFT backend, and the same image sizes.
By default, the scripts use a Gaussian PSF with radius \(r=4\) and \(\sigma=1.0\), corresponding to a kernel of size \((2r+1)\times(2r+1)=9\times 9\), together with periodic boundary conditions.
The blur strength can be adjusted by changing the PSF radius and standard deviation: increasing \(r\) enlarges the kernel support, while increasing \(\sigma\) produces a smoother and typically stronger blur.
In the scripts, these parameters can be set explicitly through the command-line arguments \texttt{--psf\_radius} and \texttt{--psf\_sigma}.
This ensures that differences in output quality and runtime reflect the restoration method rather than the surrounding implementation details.

The repository provides three levels of experimentation:
\begin{itemize}[leftmargin=1.2em]
\item \textbf{Single-image deblurring:} a direct run of the main demo script on one image, with user-specified size, regularization parameter, optional noise level, and optional PSF parameters.
\item \textbf{Lambda selection:} an auxiliary script that scans candidate values of \(\lambda\) and selects the one maximizing QSLST performance for each image/size pair.
\item \textbf{Benchmark generation:} a wrapper script that runs the full comparison across the selected Kodak images and resolutions, and saves figures, metrics, and summary tables used in the article.
\end{itemize}

A practical feature of the current demo is that it saves not only the final restored images, but also the intermediate degraded image and the blur-only image. This makes it easier to verify visually that the forward model and the restoration pipeline are consistent.

Representative commands are shown in Listing~\ref{lst:image_deblurring_demo}. An end-to-end benchmark first tunes the Tikhonov parameter \(\lambda\), then runs the full deblurring comparison using the selected values. The generated figures and JSON summaries can be used directly for quantitative analysis and paper-ready plots. For the deblurring commands below, the default blur can be overridden; for example, appending \texttt{--psf\_radius 6 --psf\_sigma 2.0} replaces the default \(9\times 9\) Gaussian PSF by a stronger blur with kernel size \(13\times 13\).

\begin{lstlisting}[language=bash,caption={Quaternion image deblurring demo and benchmark scripts.},label={lst:image_deblurring_demo}]
# From the repository root

# Single-image deblurring (no added noise)
python run_analysis.py image_deblurring \
  --image kodim16 --size 64 --lam 1e-6 \
  --ns_mode fftT --ns_iters 25 --fftT_order 2

# Single-image deblurring with additive noise
python run_analysis.py image_deblurring \
  --image kodim16 --size 64 --lam 1e-1 --snr 30 \
  --ns_mode fftT --ns_iters 25 --fftT_order 2

# Optional: tune lambda values for the benchmark images/sizes
python applications/image_deblurring/optimize_lambda.py

# Run the full benchmark
python run_analysis.py deblur_benchmark
\end{lstlisting}

\begin{figure}[t]
    \centering
    \includegraphics[width=0.8\linewidth]{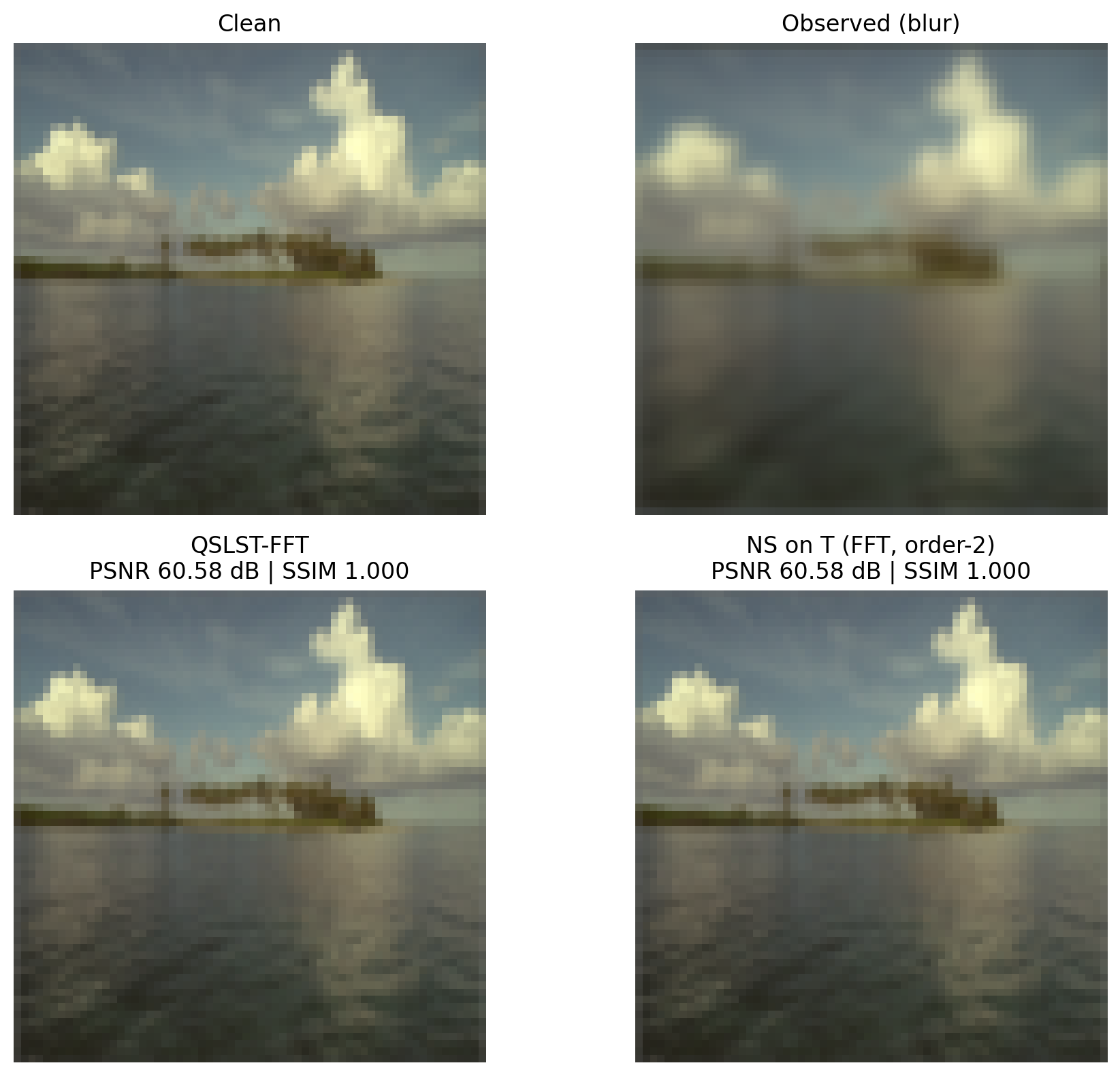}
    \caption{Example of quaternion image deblurring with \texttt{QuatIca} on \texttt{kodim16} at resolution \(64\times 64\). The grid shows the clean image, the degraded observation, and the restorations obtained with QSLST and the FFT-based Newton--Schulz method. The corresponding PSNR and SSIM values are reported in the figure title.}
    \label{fig:image_deblurring_demo}
\end{figure}

For a compact end-to-end illustration, we consider the command
\texttt{python run\_analysis.py image\_deblurring --image kodim16 --size 64 --lam 1e-6 --ns\_mode fftT --ns\_iters 25 --fftT\_order 2}.
The resulting comparison grid, shown in Figure~\ref{fig:image_deblurring_demo}, displays the clean image, the degraded observation, and the restorations obtained with QSLST and the FFT-based Newton--Schulz inverse iteration.
This example is primarily intended as a didactic illustration of the workflow and of the visual behavior of the two restoration pipelines; broader quantitative comparisons are deferred to the article~\cite{leplat2025iterativemethodscomputingmoorepenrose}.
The corresponding PSNR and SSIM values are reported directly in the figure.

\subsection{Lorenz attractor filtering as quaternion signal processing}
\label{subsec:lorenz_signal}

The Lorenz system is a classical three-dimensional nonlinear dynamical system, originally introduced by Lorenz as a simplified model for atmospheric convection \cite{lorenz1963deterministic}; see also \cite{strogatz2024nonlinear} for background on nonlinear dynamics and chaos.
It is governed by
\begin{align}
\frac{dx}{dt} &= \sigma (y - x), \\
\frac{dy}{dt} &= x (\rho - z) - y, \\
\frac{dz}{dt} &= x y - \beta z,
\end{align}
where $\sigma$, $\beta$, and $\rho$ are positive parameters. In the standard chaotic regime, one typically uses
\[
\sigma = 10, \qquad \beta = \frac{8}{3}, \qquad \rho = 28.
\]
The resulting trajectory exhibits the well-known butterfly-shaped Lorenz attractor.

To cast this problem in quaternion form, the three state components are encoded as a pure quaternion signal
\[
q(t) \;=\; x(t)\mathtt{i} + y(t)\mathtt{j} + z(t)\mathtt{k}.
\]
This representation is natural in the present framework: it packages the three coupled coordinates of the attractor into a single quaternion-valued time series while preserving their joint structure.
In the current \texttt{QuatIca} demo, a clean Lorenz trajectory is first generated numerically and then perturbed by additive Gaussian noise.
From the resulting noisy quaternion signal \(s\), the code assembles a structured quaternion linear system associated with a discrete convolution model
\[
q_i \;\approx\; \sum_{\ell=0}^{N-1} s_{i-\ell}\,h_\ell, \qquad i=0,\dots,N-1,
\]
where \(h_\ell\in\HH\) are the coefficients of an unknown right quaternion finite-impulse-response filter.
Equivalently, in matrix form, this yields a structured quaternion system
\[
A(s)h \;\approx\; q,
\]
in which \(A(s)\) is formed from shifted copies of the noisy signal.
Solving this system identifies a filter that maps the noisy trajectory to a reconstructed one, which can then be compared with the ground truth both in three-dimensional phase space and componentwise.
Additional mathematical details on the quaternion filtering viewpoint and related reconstruction methods can be found in \cite{leplat2025iterativemethodscomputingmoorepenrose,jia2021structure} and the references therein.

This example provides a compact and interpretable benchmark for quaternion structured linear solvers.
It is particularly appealing because the underlying dynamics are nonlinear and chaotic, while the reconstruction step itself is formulated as a quaternion linear algebra problem.
In the current implementation, the main reconstruction script uses the built-in Q-GMRES solver, with automatic left LU preconditioning for larger systems; see Section~\ref{subsec:QGMRES} for details.
A complementary benchmark script compares Q-GMRES with Newton--Schulz-based reconstruction and direct LU solves, thereby illustrating how different quaternion solvers behave on the same structured identification task; see also Section~\ref{subsec:NS_solver}.

\paragraph{Demonstration scripts.}
The repository provides two levels of experimentation for this application.
The command \texttt{python run\_analysis.py lorenz\_signal} runs a single reconstruction experiment and saves visualizations of the noisy and reconstructed trajectories.
Smaller values of \texttt{num\_points} are useful for debugging, whereas larger values produce smoother and more informative figures.
The command \texttt{python run\_analysis.py lorenz\_benchmark} runs a broader comparison across several problem sizes and solver choices.
In the current implementation, all generated figures are saved in the \texttt{output\_figures/} directory.

Representative commands are listed in Listing~\ref{lst:lorenz_filtering_demo}.

\begin{lstlisting}[language=bash,caption={Lorenz attractor filtering demo and benchmark.},label={lst:lorenz_filtering_demo}]
# From the repository root
# Fast Lorenz filtering demo
python run_analysis.py lorenz_signal --num_points 100

# Default / higher-quality Lorenz filtering demo
python run_analysis.py lorenz_signal
python run_analysis.py lorenz_signal --num_points 500

# Benchmark comparison of several solvers
python run_analysis.py lorenz_benchmark

# Newton--Schulz-only benchmark on a small test case
python run_analysis.py lorenz_benchmark --methods newton --points 100
\end{lstlisting}

\begin{figure}[ht!]
    \centering

    \begin{minipage}{0.48\linewidth}
        \centering
        \includegraphics[width=\linewidth]{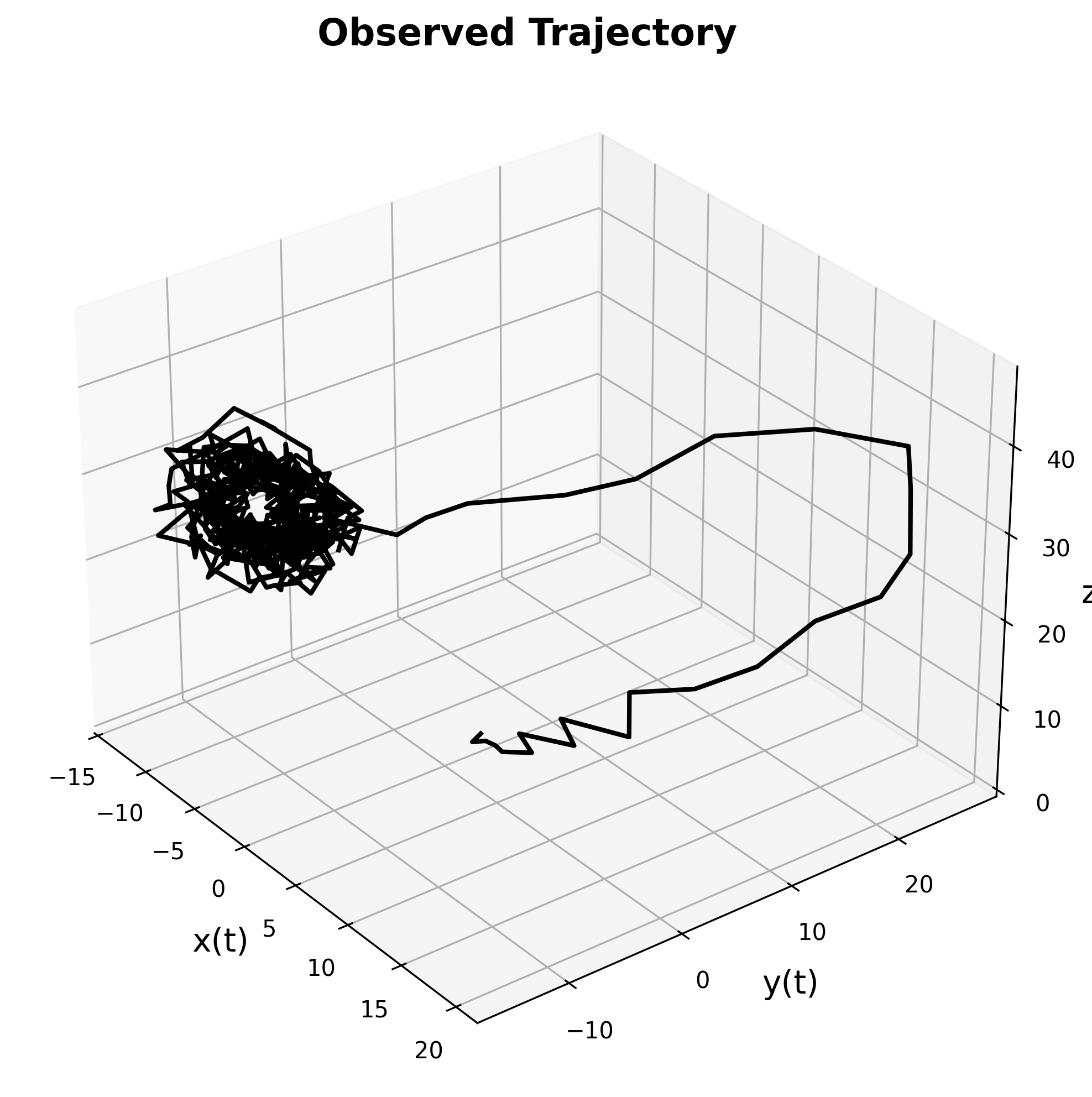}
    \end{minipage}\hfill
    \begin{minipage}{0.48\linewidth}
        \centering
        \includegraphics[width=\linewidth]{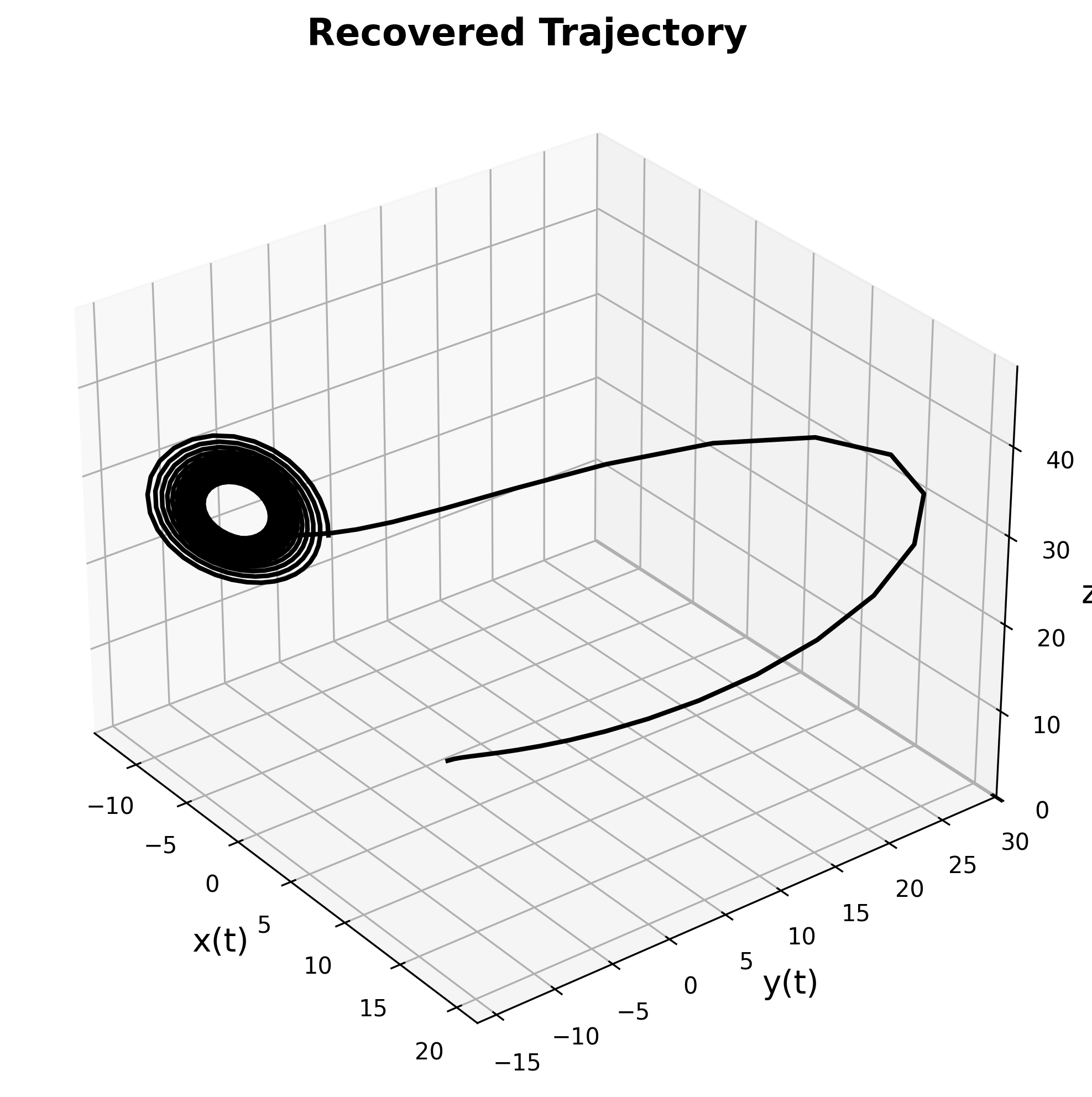}
    \end{minipage}

    \vspace{0.5em}

    \begin{minipage}{0.48\linewidth}
        \centering
        \includegraphics[width=\linewidth]{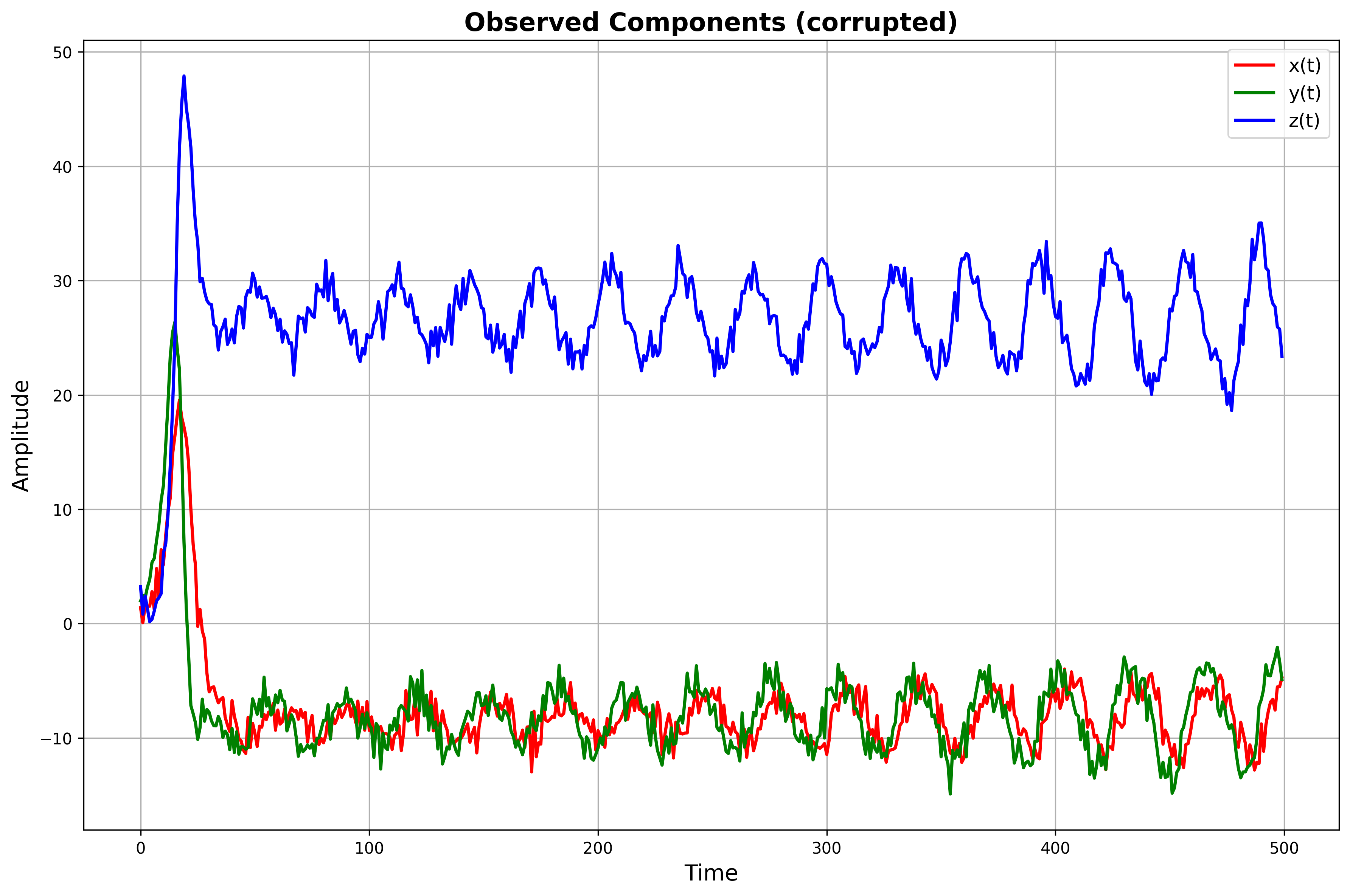}
    \end{minipage}\hfill
    \begin{minipage}{0.48\linewidth}
        \centering
        \includegraphics[width=\linewidth]{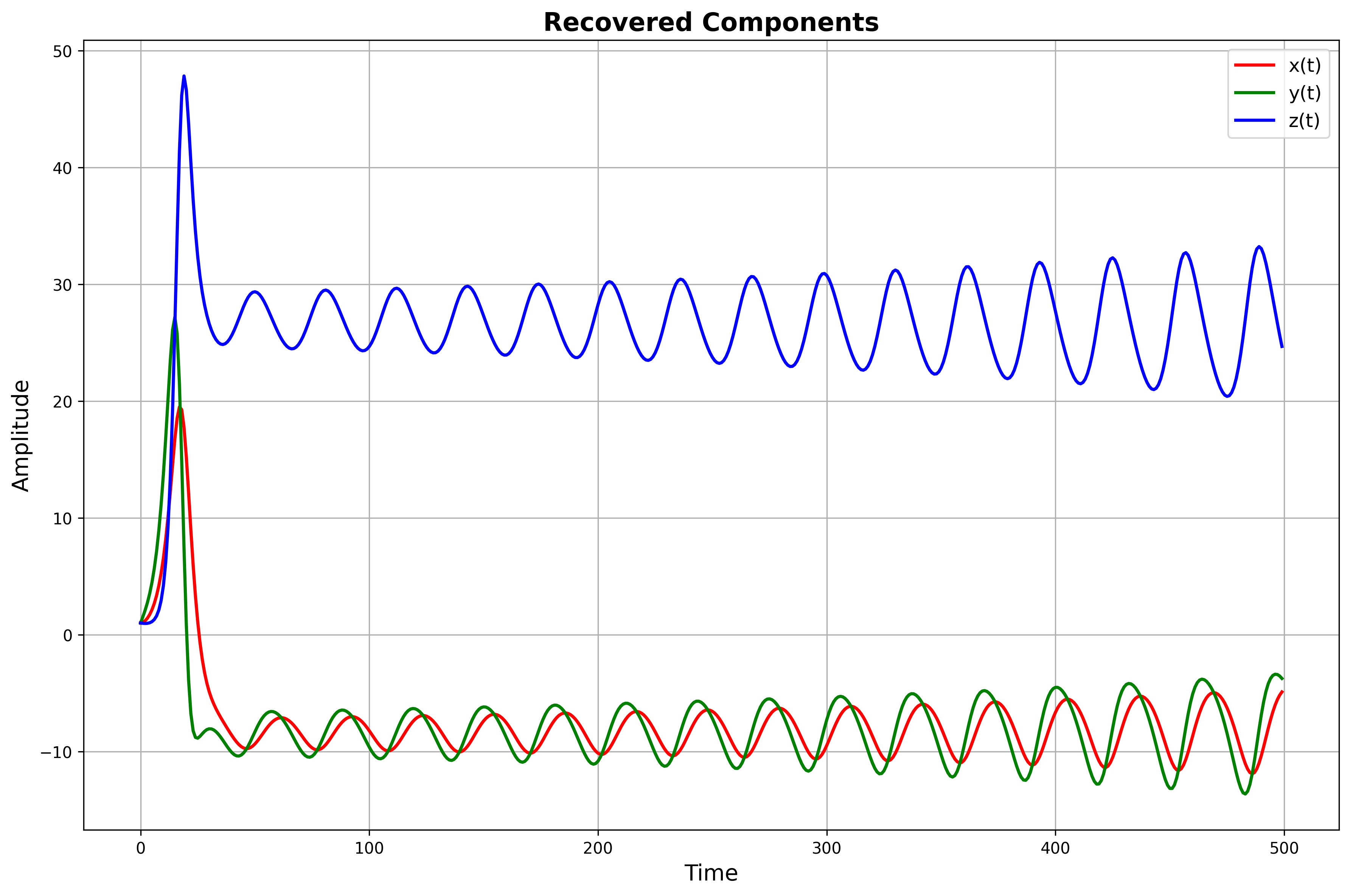}
    \end{minipage}

    \caption{Lorenz attractor filtering example generated with \texttt{python run\_analysis.py lorenz\_signal --num\_points 500}. Top: noisy observed trajectory (left) and reconstructed trajectory (right) in three-dimensional phase space. Bottom: observed component signals (left) and reconstructed component signals (right). The reconstruction is obtained with the built-in Q-GMRES solver with automatic left LU preconditioning; see Section~\ref{subsec:QGMRES}. A complementary Newton--Schulz-based experiment can be generated with \texttt{python run\_analysis.py lorenz\_benchmark --methods newton --points 500}; see Section~\ref{subsec:NS_solver}.}
    \label{fig:lorenz_signal_demo}
\end{figure}

For a compact end-to-end illustration, we consider the command
\texttt{python run\_analysis.py lorenz\_signal --num\_points 500}.
The corresponding outputs are shown in Figure~\ref{fig:lorenz_signal_demo}.
The top row displays the noisy observed trajectory and the reconstructed three-dimensional trajectory, while the bottom row shows the corresponding component signals.
In this example, the reconstruction is obtained with the built-in Q-GMRES solver with automatic left LU preconditioning.
This figure is primarily intended as a didactic illustration of the workflow and of the visual quality of the recovered Lorenz signal; additional solver comparisons, including Newton--Schulz-based reconstructions, can be generated with the benchmark command above and are discussed further in \cite{leplat2025iterativemethodscomputingmoorepenrose}.

\newpage
\subsection{Quaternion Trajectory Interpolation and Analysis}
\label{subsec:qtraj_demo}

To illustrate the \texttt{qtraj} module, \texttt{QuatIca} includes a reproducible demo comparing three
practical interpolation schemes for quaternion-valued orientation keyframes on \(S^3\):
piecewise SLERP, SQUAD, and log--exp spline interpolation.
The goal of this experiment is not only to visualize the resulting trajectories, but also to compare
their smoothness through simple kinematic diagnostics derived from discrete angular velocity.

The demo first generates a small sequence of unit-quaternion keyframes
\[
q_0,q_1,\dots,q_{K-1}\in S^3
\]
at uniformly spaced times
\[
0=t_0<t_1<\cdots<t_{K-1}=1.
\]
The keyframes are built from a smooth synthetic axis--angle construction and then processed with
\texttt{enforce\_sign\_continuity} to remove artificial sign flips due to the double-cover relation
\(q\sim -q\).
Three interpolants are then sampled on a dense time grid:
\begin{itemize}[leftmargin=1.2em]
\item \emph{piecewise SLERP}, used as a geodesic baseline on each interval \([t_i,t_{i+1}]\);
\item \emph{SQUAD}, which introduces control quaternions to obtain smoother transitions between segments;
\item \emph{log--exp cubic spline}, which maps the keyframes to \(\R^3\) via the quaternion logarithm,
fits a cubic spline there, and maps the result back to \(S^3\) with the quaternion exponential.
\end{itemize}

For each sampled trajectory, the demo estimates a discrete angular velocity using the midpoint log-map formula
\[
\omega(t_{k+\frac12}) \approx \frac{2}{\Delta t}\,\log\!\big(q(t_{k+1})\,q(t_k)^{-1}\big),
\qquad
t_{k+\frac12}=\frac{t_k+t_{k+1}}{2},
\]
and reports several diagnostic quantities:
a heuristic smoothness energy based on discrete angular acceleration,
the root-mean-square and maximal angular speed,
and, for piecewise constructions, the estimated jump of angular velocity at keyframes.
These quantities provide a simple quantitative comparison of how smoothly the different interpolants
connect the same orientations.

\paragraph{Demonstration script.}
The quaternion interpolation demo can be run directly from the repository.
A typical command is:

\begin{lstlisting}[language=bash,caption={Quaternion geodesic interpolation demo.},label={lst:qtraj_demo}]
# From the repository root
python run_analysis.py qtraj_demo
\end{lstlisting}
The exact script path can be adjusted to the local repository layout.

\paragraph{Produced outputs.}
The demo saves a set of figures in the output directory associated with quaternion geodesic interpolation:
\begin{itemize}[leftmargin=1.2em]
\item \texttt{angular\_speed.(png/pdf)}: profile of \(\|\omega(t)\|\) for the three interpolants;
\item \texttt{angular\_accel.(png/pdf)}: discrete angular-acceleration magnitude
\(\|\dot\omega(t)\|\);
\item \texttt{keyframe\_errors.(png/pdf)}: geodesic interpolation error at the prescribed keyframe times;
\item \texttt{s2\_trajectory.(png/pdf)}: trajectory on \(S^2\) obtained by rotating a fixed unit vector
by the interpolated orientations.
\end{itemize}

\paragraph{What the demo shows.}
The angular-speed and angular-acceleration plots highlight the main qualitative difference between the methods.
Piecewise SLERP is geodesically natural on each segment, but typically exhibits visible changes in angular velocity
at the knot points because it is only \(C^0\) as a multi-segment trajectory.
By contrast, SQUAD and the log--exp spline usually reduce these visible discontinuities on the smooth synthetic keyframes used in the default demo.
In the current controlled experiments, their angular-speed profiles are broadly comparable, with SQUAD retaining a local spline construction and the log--exp spline often appearing slightly smoother at a global level.
The latter, however, remains more dependent on the choice of logarithmic chart.

The keyframe-error plot serves as a consistency check: all three methods interpolate the prescribed keyframes
up to numerical precision.
Finally, the \(S^2\) visualization provides an intuitive geometric picture of the motion by showing how a fixed
unit vector is transported by the interpolated orientations.
Even when all methods pass through the same keyframes, their between-keyframe paths on \(S^2\) can differ
significantly, making the effect of the interpolation model easy to interpret visually.

\begin{figure}[ht!]
  \centering
  \includegraphics[width=0.77\linewidth]{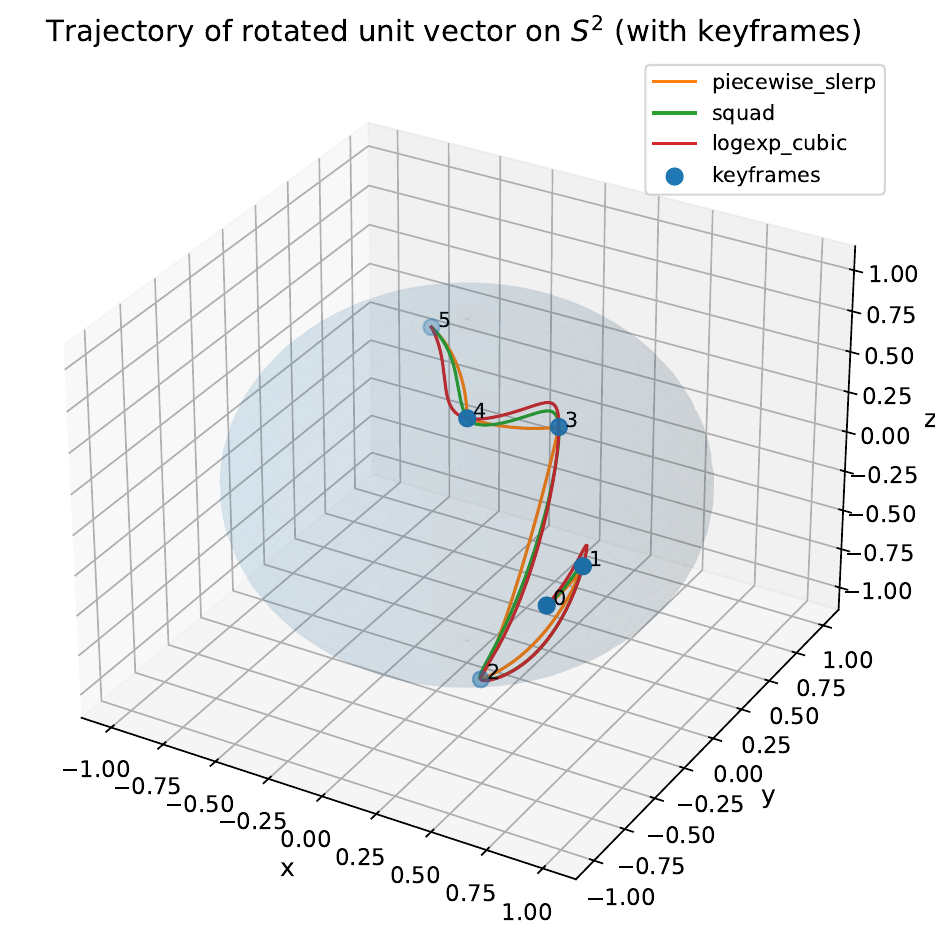}
  \caption{Trajectory on \(S^2\) obtained by rotating a fixed unit vector with the interpolated quaternion orientations. Using the same randomly generated keyframes, piecewise SLERP, SQUAD, and log--exp spline interpolation all produce valid rotation trajectories passing through the prescribed keyframes, while following different between-keyframe paths on the sphere.}
  \label{fig:qtraj_s2_trajectory}
\end{figure}

As a representative visual output, we include in Figure~\ref{fig:qtraj_s2_trajectory}
the trajectory on \(S^2\) obtained by rotating a fixed unit vector by the interpolated
orientations. For a smooth synthetic set of keyframes, all three methods produce valid
rotation trajectories while exhibiting clearly different between-keyframe paths.

\section{Testing, numerical validation, and benchmarks}
\label{sec:repro}
\subsection{Testing and numerical validation}
\label{subsec:testing_validation}

A substantial part of \texttt{QuatIca}'s development effort is devoted to testing and numerical validation.
The library currently includes an extensive collection of unit tests covering the main routines and helper
functions implemented across the package. These tests are designed to verify both algebraic correctness and
numerical behavior, including Hermitian preservation, decomposition residuals, orthogonality/unitarity checks,
solver convergence, and consistency of quaternion-specific matrix operations. In the repository, the test files
follow a structured naming convention of the form \texttt{test\_<function\_name>\_<details>.py}, and can be run
directly with the \texttt{pytest} framework to obtain detailed reports.

Since the test suite currently contains a large number of unit and regression tests, it is neither practical nor
particularly informative to report them exhaustively in the manuscript. Instead, in this section we highlight the
most representative advanced numerical validation currently available in \texttt{QuatIca}, namely the end-to-end
validation of the \texttt{OptiQ} module. This example is especially informative because it simultaneously tests
several critical layers of the implementation: quaternion Hermitian operators, hat-space conditioning, Newton
steps, barrier continuation, and recovery of a certified known optimum. As such, it provides a compact but
stringent illustration of the numerical reliability targeted throughout the library.

\subsubsection{Validation of \texttt{OptiQ}: certified known-optimum, trace-bounded regression instance}
\label{subsec:optiq_known_optimum_demo}
To validate OptiQ end-to-end (operators, hat-space conditioning, Newton step, and $\mu$-continuation),
we build an SDP instance with a certified known optimum $X_\star$ and an explicit dual certificate.
This construction is used in the script \texttt{tests/optiQ/run\_optiq\_known\_optimum.py}
and produces reproducible figures for the article.

\paragraph{Step 1: complementary primal/dual projectors.}
Choose a rank $r\in\{1,2,\dots,n-1\}$ and a quaternion unitary matrix $V$.
Define the rank-$r$ projector
\[
X_\star := V
\begin{bmatrix}
I_r & 0\\[2pt]
0 & 0
\end{bmatrix}
V^H,
\]
and a complementary PSD matrix
\[
C := V
\begin{bmatrix}
0 & 0\\[2pt]
0 & I_{n-r}
\end{bmatrix}
V^H.
\]
Then $X_\star\succeq 0$, $C\succeq 0$, and $\ip{C}{X_\star}=0$. Moreover, for any $X\succeq 0$,
$\ip{C}{X}\ge 0$, so the optimal value is $0$ whenever $X_\star$ is feasible.

\paragraph{Step 2: a trace constraint to prevent blow-up.}
Barrier methods can exhibit eigenvalue inflation if the feasible set contains ``free'' PSD directions
that do not increase the linear objective. To bound the feasible set, we impose
\[
\tr(X)=\tr(X_\star)=r,
\]
implemented as a linear constraint $\ip{H_0}{X}=b_0$ with
\[
H_0 := \frac{1}{\sqrt{n}}I,\qquad b_0 := \ip{H_0}{X_\star}=\frac{r}{\sqrt{n}}.
\]
The normalization $1/\sqrt{n}$ is used for conditioning.

\paragraph{Step 3 (optional): additional trace-free constraints.}
To reduce degeneracy while preserving boundedness, we add $m_{\rm extra}$ random trace-free constraints:
sample Hermitian $H$, then remove its trace component by
\[
H \leftarrow \herm\!\left(H - \alpha I\right),
\qquad
\alpha := \frac{\ip{H}{I}}{\ip{I}{I}}.
\]
We then set $b_i:=\ip{H_i}{X_\star}$ so that $A(X_\star)=b$ holds by construction.
Finally, we build hat operators via Gram--Cholesky and define $\hat b := \hat A(X_\star)$.

\paragraph{Step 4: strictly feasible warm start (Slater point).}
A simple strictly feasible (positive definite) warm start preserving the trace is
\[
X_{\mathrm{feas}} := (1-\varepsilon)\,X_\star + \varepsilon\,\frac{r}{n}I \ \succ\ 0,
\qquad \tr(X_{\mathrm{feas}})=r,
\]
with $\varepsilon\in(0,1)$ small (e.g., $\varepsilon=10^{-2}$).

\paragraph{What we expect to observe.}
Along the barrier path $X(\mu)$, we typically observe:
(i) hat-space feasibility $\|\hat b-\hat A(X(\mu))\|_2$ near machine precision,
(ii) objective scaling $\ip{C}{X(\mu)}\approx (n-r)\mu$ for small $\mu$,
and (iii) spectral separation with $\lambda_{\min}(X(\mu))\asymp \mu$.

\paragraph{Figures produced by the demo script.}

The demo script can be run as described in Listing~\ref{lst:optiq_demo}.
\begin{lstlisting}[caption={Running the certified known-optimum OptiQ regression demo.},label={lst:optiq_demo}]
# from the repo root
python tests/optiQ/run_optiq_known_optimum.py
# or equivalently from the wrapper run_analysis
python run_analysis.py optiq_known
# figures saved under validation_output/optiQ/known_optimum_trace/
\end{lstlisting}

The demo script saves, for a sequence of decreasing $\mu$ values:
\begin{itemize}[leftmargin=1.2em]
\item \texttt{obj\_vs\_mu.png}: $\ip{C}{X(\mu)}$ vs.\ $\mu$ (log--log, $\mu$ axis inverted);
\item \texttt{relerr\_vs\_mu.png}: $\|X(\mu)-X_\star\|_F/\|X_\star\|_F$ vs.\ $\mu$;
\item \texttt{eigs\_vs\_mu.png}: $\lambda_{\min}(X(\mu))$ and $\lambda_{\max}(X(\mu))$ vs.\ $\mu$;
\item \texttt{abs\_heatmaps\_full.png}: heatmaps of $|X_\star|$, $|X(\mu_{\min})|$, and $|X_\star-X(\mu_{\min})|$;
\item \texttt{block\_$k$x$k$\_components.png}: side-by-side $k\times k$ blocks of the real/$\mathtt{i}$/$\mathtt{j}$/$\mathtt{k}$
components of $X_\star$ and $X(\mu_{\min})$ (and their difference), for a small block size $k$ (default $k=5$).
\end{itemize}

In this article, we include a small subset of these outputs as a compact end-to-end validation of \texttt{OptiQ}.
Together, the plots confirm the expected barrier-path behavior: as the barrier parameter $\mu$ decreases,
the computed solution $X(\mu)$ approaches the certified optimum $X_\star$ while feasibility is maintained.

\begin{figure}[ht!]
\centering

\begin{subfigure}{\linewidth}
  \centering
  \includegraphics[width=0.99\linewidth]{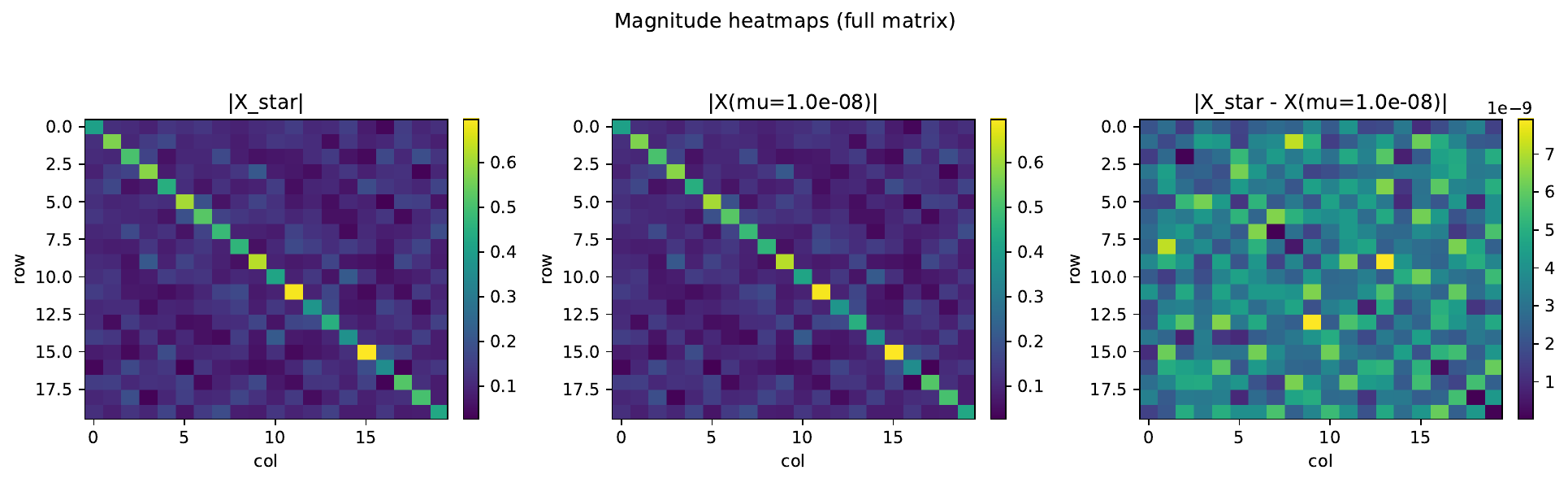}
  \caption{Absolute-value heatmaps comparing $|X_\star|$, $|X(\mu_{\min})|$, and $|X_\star-X(\mu_{\min})|$.}
  \label{fig:optiq_abs_heatmaps_full}
\end{subfigure}

\vspace{0.6em}

\begin{subfigure}{\linewidth}
  \centering
  \includegraphics[width=0.60\linewidth]{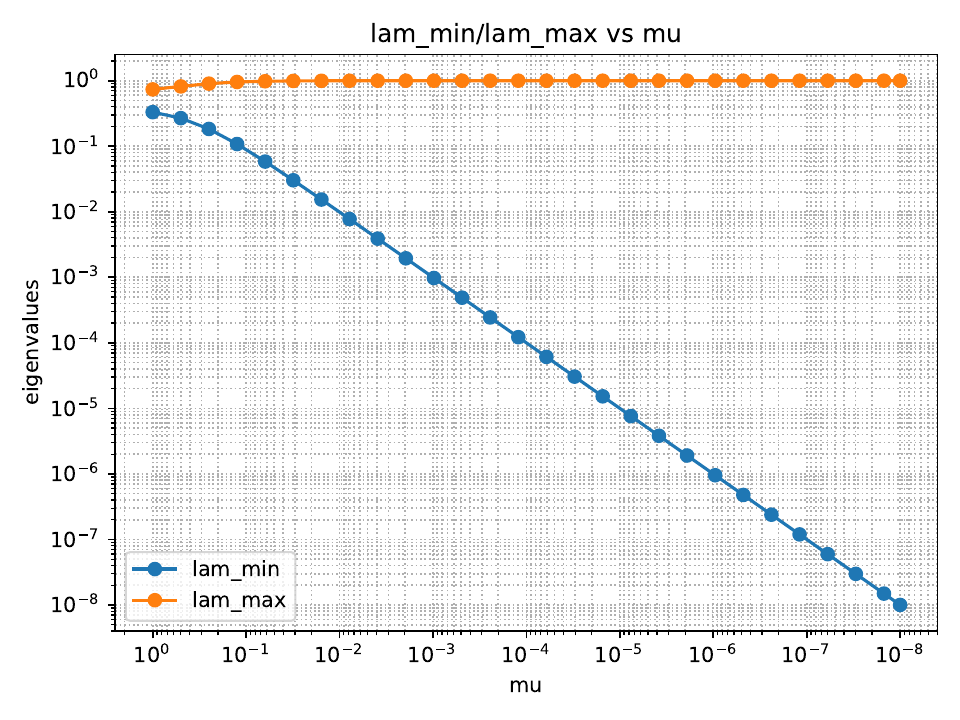}
  \caption{Extreme eigenvalues along the barrier path: $\lambda_{\min}(X(\mu))\to 0$ and $\lambda_{\max}(X(\mu))\to 1$.}
  \label{fig:optiq_eigs_vs_mu}
\end{subfigure}

\vspace{0.6em}

\begin{subfigure}{0.49\linewidth}
  \centering
  \includegraphics[width=\linewidth]{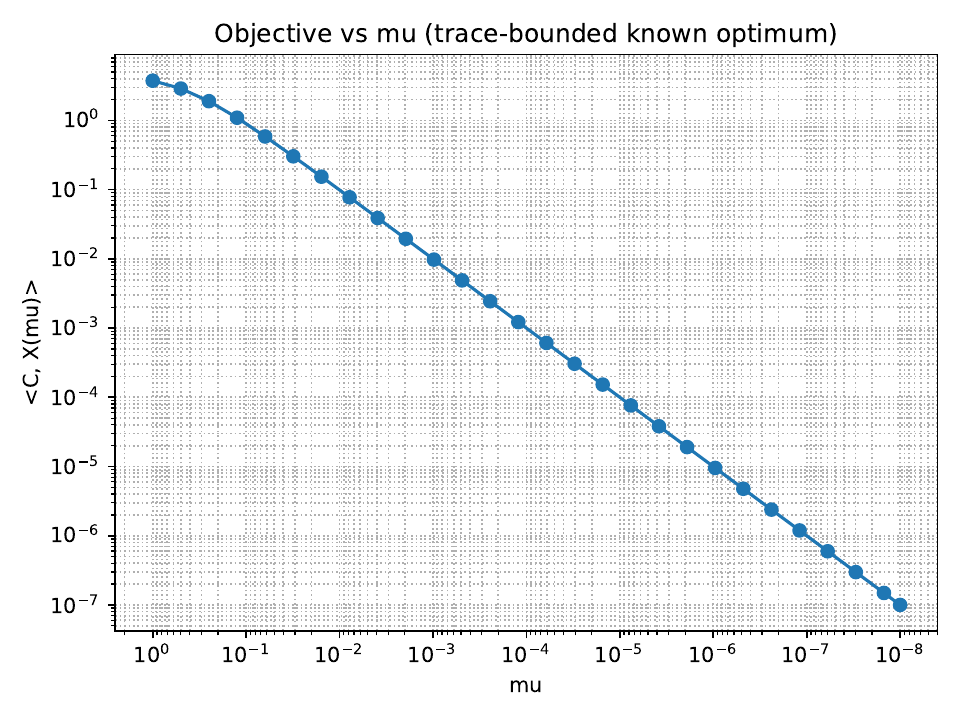}
  \caption{Objective value $\ip{C}{X(\mu)}$ versus $\mu$ (log--log), showing the expected decrease as $\mu$ is reduced.}
  \label{fig:optiq_obj_vs_mu}
\end{subfigure}\hfill
\begin{subfigure}{0.49\linewidth}
  \centering
  \includegraphics[width=\linewidth]{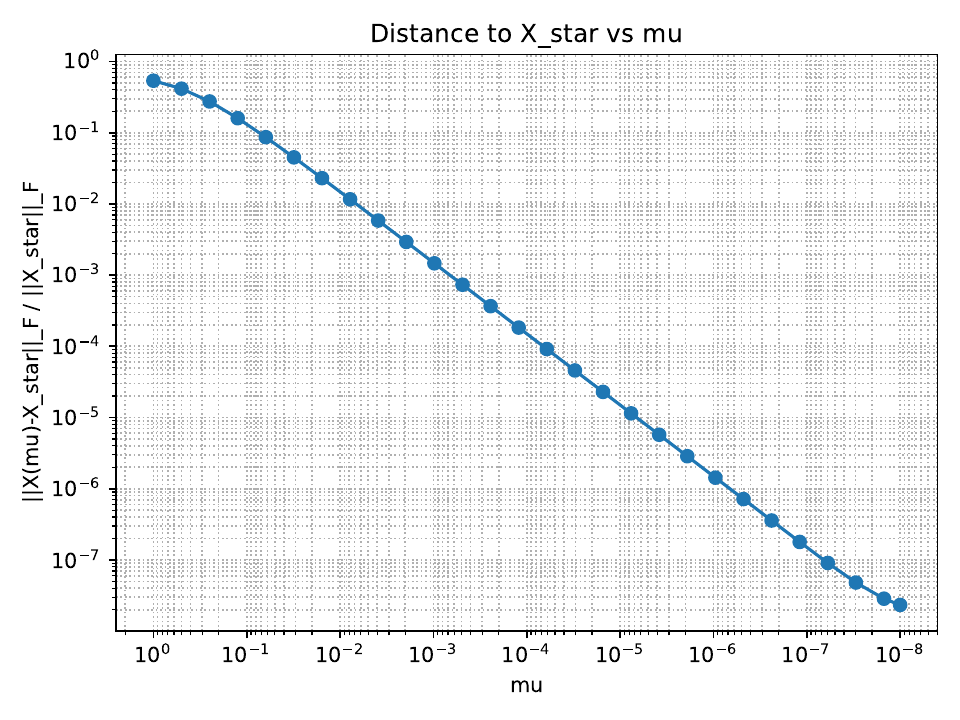}
  \caption{Relative Frobenius error $\|X(\mu)-X_\star\|_F/\|X_\star\|_F$ versus $\mu$, tending to zero until the stopping threshold.}
  \label{fig:optiq_relerr_vs_mu}
\end{subfigure}

\caption{Certified known-optimum OptiQ regression instance: representative diagnostics along the barrier path.}
\label{fig:optiq_known_optimum_diagnostics}
\end{figure}

Concretely, as one can observe in Figures~\ref{fig:optiq_abs_heatmaps_full}--\ref{fig:optiq_eigs_vs_mu},
the diagnostics provide complementary checks of correctness. The heatmaps in
Figure~\ref{fig:optiq_abs_heatmaps_full} provide a direct visual comparison showing that $|X(\mu_{\min})|$
matches $|X_\star|$ for the smallest value of $\mu$ considered. The curves in
Figures~\ref{fig:optiq_obj_vs_mu} and \ref{fig:optiq_relerr_vs_mu} show that the objective value
$\ip{C}{X(\mu)}$ and the relative Frobenius error $\|X(\mu)-X_\star\|_F/\|X_\star\|_F$ decrease with $\mu$
until the stopping threshold (typically $10^{-8}$). Finally, Figure~\ref{fig:optiq_eigs_vs_mu} tracks the
extreme eigenvalues of $X(\mu)$ and exhibits the expected spectral trend: $\lambda_{\min}(X(\mu))$ decreases
toward $0$ as $\mu\downarrow 0$, while $\lambda_{\max}(X(\mu))$ approaches $1$ (consistent with the
projector structure of $X_\star$). Overall, these results provide clear evidence that the implementation
reproduces the known optimum and follows the intended continuation trajectory.

\subsection{Benchmarks}
\label{subsec:benchmarks}

Beyond algebraic validation, unit tests and end-to-end application demos, \texttt{QuatIca} provides dedicated benchmark
experiments to assess the practical performance of its main numerical solvers. These tests are designed to show
how the methods behave as the problem dimension increases and as the underlying linear systems or matrices become
more challenging. Depending on the routine, we report quantities such as runtime, iteration count, approximation
accuracy, residual decay, and robustness across representative scenarios.

The benchmark suite considered in this article focuses on three important families of routines:
quaternion low-rank and spectral decompositions,
Newton--Schulz pseudoinverse algorithms,
and Q-GMRES with and without built-in preconditioning.
Each benchmark is presented separately below together with a short description of the test protocol and a subset
of representative figures generated by the corresponding analysis scripts.

\subsubsection{Q-SVD benchmark: low-rank approximation accuracy and efficiency}
\label{subsec:qsvd_benchmark}

To assess the practical behavior of the quaternion SVD routines implemented in \texttt{QuatIca},
we compare three representative variants:
the deterministic truncated Q-SVD, a randomized Q-SVD based on Gaussian sketching and power iterations,
and the pass-efficient Q-SVD variant.
The benchmark is organized around three complementary tests:
(i) an exact low-rank recovery sanity check,
(ii) approximation error as a function of the target rank,
and (iii) runtime scalability as the matrix size increases. 

In the first experiment, the script constructs an exact rank-10 quaternion matrix
\[
X = AB,\qquad A\in\HH^{100\times 10},\quad B\in\HH^{10\times 200},
\]
and compares the three methods when the requested target rank is \(R=10\).
In the second experiment, it considers a rectangular matrix with smooth spectral decay,
of size \(300\times 200\), and evaluates the relative reconstruction error
\[
\frac{\|X-\hat X_R\|_F}{\|X\|_F},
\qquad
\hat X_R = U_R \Sigma_R V_R^H,
\]
for a range of target ranks \(R\in\{5,10,20,40,60,80,100\}\).
Finally, the runtime experiment measures the wall-clock cost of the three methods on increasingly large dense quaternion matrices, with default sizes ranging from \(100\times 80\) up to \(1200\times 900\), while keeping the target rank fixed. The benchmark exports a full dashboard together with separate PDF panels and a CSV file of aggregated metrics. 

\paragraph{Demonstration script.}
The benchmark can be launched directly from the analysis wrapper:
\begin{lstlisting}[language=bash,caption={Q-SVD benchmark for deterministic and randomized variants.},label={lst:qsvd_bench}]
# from the repo root
python run_analysis.py qsvd_bench

# Example: increase the pass budget of the pass-efficient variant
python run_analysis.py qsvd_bench --n-passes 4

# Show all available options
python run_analysis.py qsvd_bench --help
\end{lstlisting}

For more aggressive runtime-oriented stress tests, the randomized parameters can also be overridden from the command line.
In particular, the pass budget of the pass-efficient variant can be controlled by appending
\texttt{--n-passes N\_PASSES},
which typically improves its approximation quality and can make its accuracy nearly match that of the standard randomized Q-SVD.
Additional options can be inspected with \texttt{--help}.

\begin{figure}[ht!]
  \centering

  \begin{subfigure}{0.45\linewidth}
    \centering
    \includegraphics[width=\linewidth]{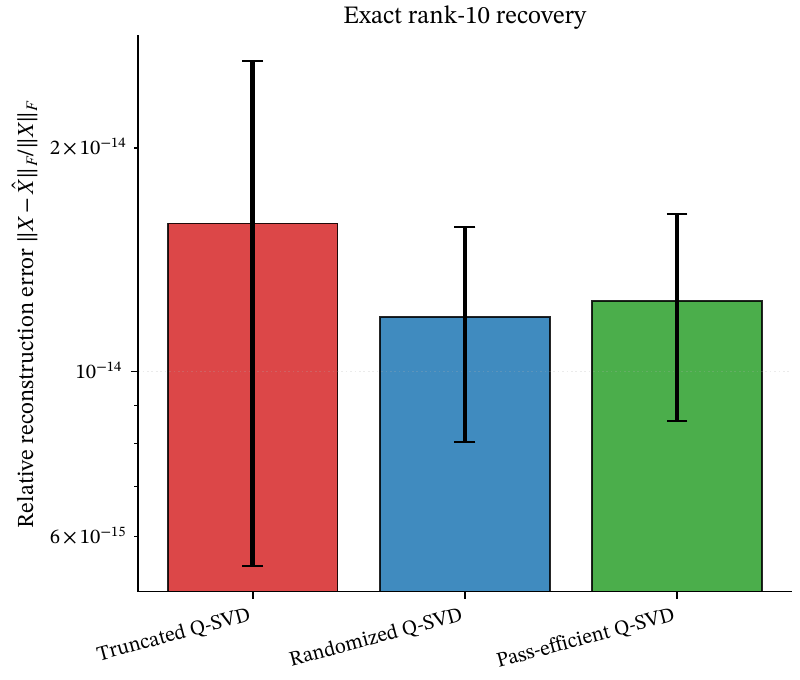}
    \caption{Exact rank-10 recovery.}
    \label{fig:qsvd_exact_rank10}
  \end{subfigure}\hfill
  \begin{subfigure}{0.45\linewidth}
    \centering
    \includegraphics[width=\linewidth]{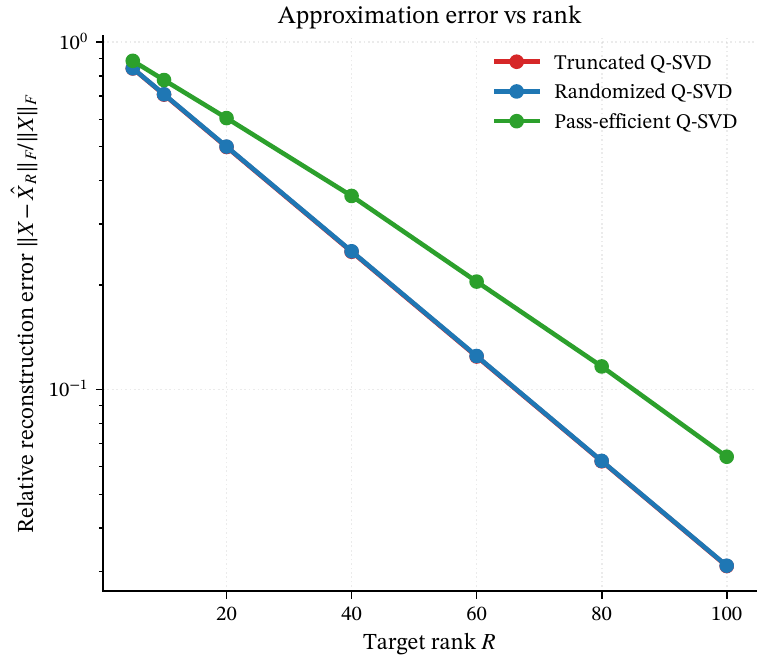}
    \caption{Approx. error vs target rank.}
    \label{fig:qsvd_error_vs_rank}
  \end{subfigure}\hfill
  \begin{subfigure}{0.45\linewidth}
    \centering
    \includegraphics[width=\linewidth]{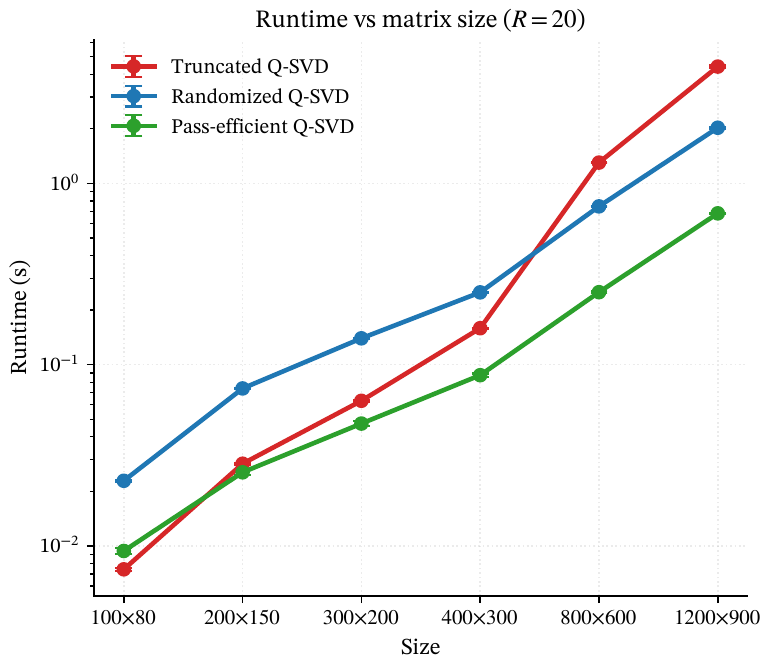}
    \caption{Runtime versus matrix size.}
    \label{fig:qsvd_runtime_vs_size}
  \end{subfigure}

  \caption{Representative Q-SVD benchmark outputs generated by the paper-ready benchmark script. The top left panel verifies near machine-precision recovery on an exact rank-10 quaternion matrix when the requested rank is correct. The top right panel shows that the relative reconstruction error decreases as the target rank increases and that the randomized variants closely follow the deterministic truncated Q-SVD. The bottom panel illustrates the increasing runtime advantage of the randomized approaches as the matrix size grows, with the pass-efficient variant becoming especially competitive on the largest tested problems.}
  \label{fig:qsvd_benchmark}
\end{figure}

\paragraph{What the benchmark shows.}
A representative subset of the outputs is displayed in Figure~\ref{fig:qsvd_benchmark}.
Figure~\ref{fig:qsvd_exact_rank10} shows that all three methods recover an exact rank-10 quaternion matrix to near machine precision when the target rank is chosen correctly.
Figure~\ref{fig:qsvd_error_vs_rank} then illustrates the expected decrease of the relative reconstruction error as the target rank increases, with both randomized methods remaining close to the deterministic truncated Q-SVD.
When the pass-efficient variant is slightly less accurate, this difference is mainly due to its more restrictive pass budget: in the current benchmark configuration, it uses only two passes over the data by default.
This choice favors efficiency, but may mildly reduce approximation quality.
In practice, the accuracy gap can be reduced significantly by allowing a slightly larger number of passes, bringing the pass-efficient variant much closer to the Gaussian-sketching version.
Finally, Figure~\ref{fig:qsvd_runtime_vs_size} highlights the practical interest of randomized low-rank approximation in the quaternion setting: as the matrix size grows, the runtime of the full truncated Q-SVD increases more rapidly, while the randomized and pass-efficient variants remain substantially cheaper.
In particular, the pass-efficient method becomes the most attractive option among the tested implementations for the largest matrix sizes considered here.

\subsubsection{Newton--Schulz MP pseudoinverse benchmark}
\label{subsec:ns_benchmark}

\texttt{QuatIca} implements several Newton--Schulz-type iterations for approximating the MP
pseudoinverse of a quaternion matrix using only matrix--matrix multiplications. These include the classical
Newton--Schulz iteration, a damped variant, and a higher-order third-order scheme. A broader numerical study of
these methods and their variants is presented in our companion work~\cite{leplat2025iterativemethodscomputingmoorepenrose}.
Here, we restrict attention to a compact benchmark illustrating the practical behavior of three representative
variants available in \texttt{QuatIca}:
\begin{itemize}[leftmargin=1.2em]
\item classical Newton--Schulz, denoted \(\mathrm{NS}\;(\gamma=1)\);
\item a damped Newton--Schulz variant, denoted \(\mathrm{NS}\;(\gamma=\tfrac12)\);
\item a third-order higher-order Newton--Schulz scheme.
\end{itemize}

The benchmark considers two matrix families:
(i) dense random quaternion matrices, and
(ii) rectangular ill-conditioned matrices built with controlled spectral decay.
Several matrix sizes are tested, and each configuration is repeated across multiple random seeds in order to
obtain aggregated runtime and convergence statistics.

For each run, the script records the total runtime, the number of performed iterations, and the four standard
Penrose-type diagnostics
\[
E_1=\|AXA-A\|_F,\qquad
E_2=\|XAX-X\|_F,\qquad
E_3=\|(AX)^H-AX\|_F,\qquad
E_4=\|(XA)^H-XA\|_F,
\]
where \(X\) denotes the computed pseudoinverse approximation.
In addition, the benchmark stores the full convergence history of the primary residual
\[
E_1^{(k)}=\|A X_k A-A\|_F,
\]
which is used to compare the iteration profiles of the three methods.
A run is declared successful when the final value of \(E_1\) reaches the prescribed target tolerance.

\paragraph{Demonstration script.}
The full benchmark can be launched directly from the analysis wrapper:
\begin{lstlisting}[language=bash,caption={Newton--Schulz pseudoinverse benchmark.},label={lst:ns_bench}]
# from the repo root
python run_analysis.py ns_compare
\end{lstlisting}

The script exports a dashboard figure together with separate PDF panels for article inclusion, as well as raw CSV/NPZ
data for reproducibility.

\begin{figure}[ht!]
  \centering

  \begin{subfigure}{0.49\linewidth}
    \centering
    \includegraphics[width=\linewidth]{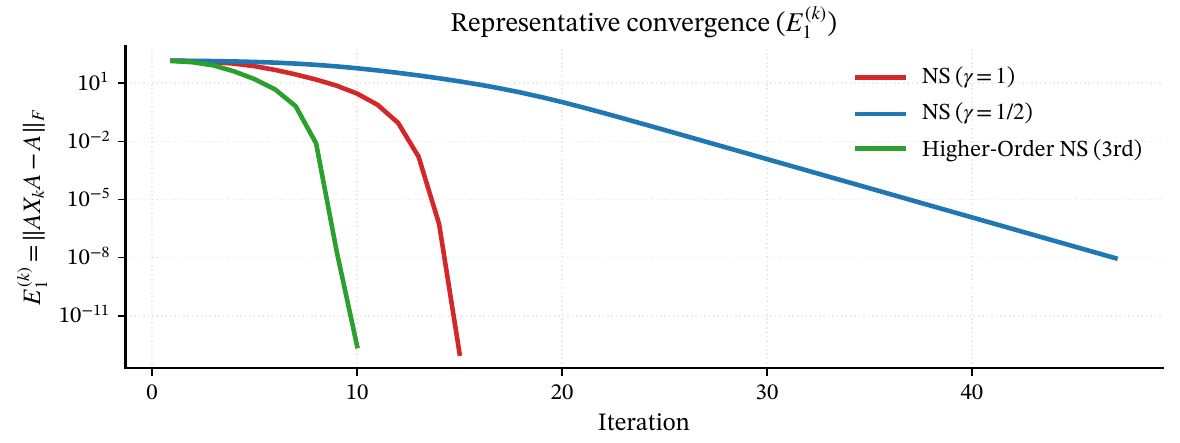}
    \caption{Representative convergence history of \(E_1^{(k)}=\|AX_kA-A\|_F\).}
    \label{fig:ns_convergence_example}
  \end{subfigure}\hfill
  \begin{subfigure}{0.49\linewidth}
    \centering
    \includegraphics[width=\linewidth]{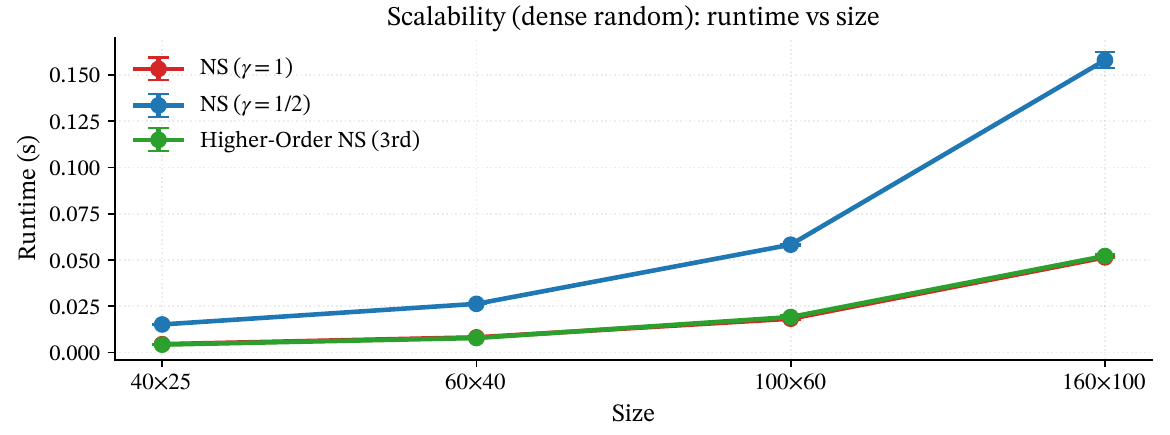}
    \caption{Average runtime as a function of matrix size.}
    \label{fig:ns_time_vs_size}
  \end{subfigure}

  \vspace{0.6em}

  \begin{subfigure}{0.58\linewidth}
    \centering
    \includegraphics[width=\linewidth]{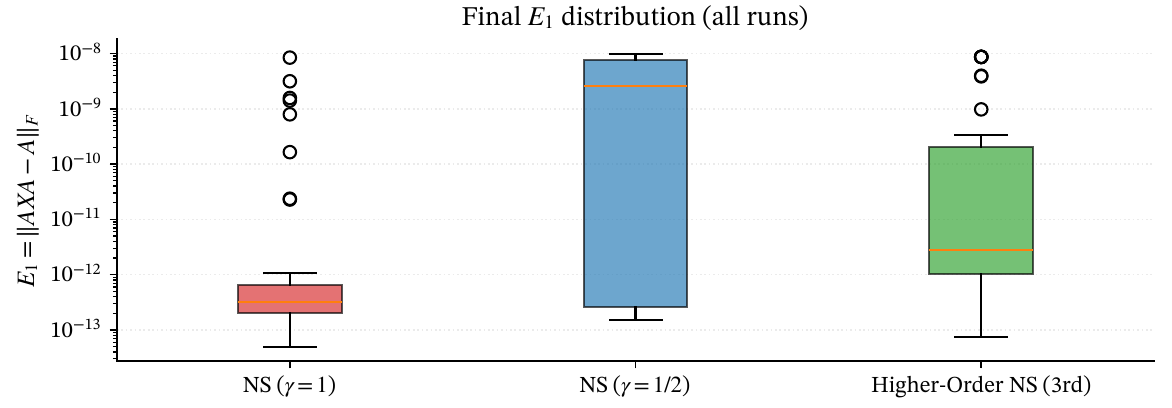}
    \caption{Distribution of the final residual \(E_1\) across all runs.}
    \label{fig:ns_final_residual_boxplot}
  \end{subfigure}

  \caption{Representative benchmark outputs for the Newton--Schulz pseudoinverse solvers implemented in \texttt{QuatIca}. The benchmark compares the classical Newton--Schulz iteration, a damped variant, and a third-order higher-order scheme on dense random and ill-conditioned quaternion matrices of increasing size.}
  \label{fig:ns_benchmark}
\end{figure}

\paragraph{What the benchmark shows.}
A representative subset of the outputs is displayed in Figure~\ref{fig:ns_benchmark}.
Figure~\ref{fig:ns_convergence_example} compares the convergence histories of the three variants on a representative
dense random matrix through the residual \(E_1^{(k)}\), while Figure~\ref{fig:ns_time_vs_size} reports the average
runtime as the matrix size increases. Figure~\ref{fig:ns_final_residual_boxplot} complements these plots by
summarizing the final residual distribution across all runs. Together, these results provide a compact overview of
the practical behavior of the Newton--Schulz pseudoinverse variants currently available in \texttt{QuatIca}.

\subsubsection{Q-GMRES benchmark: impact of built-in LU preconditioning}
\label{subsec:qgmres_benchmark}

To assess the robustness and scalability of the quaternion GMRES implementation in \texttt{QuatIca},
we benchmark Q-GMRES on synthetic linear systems
\[
Ax=b,\qquad b=Ax_\star,
\]
where \(x_\star\) is a known reference solution.
The benchmark compares two solver configurations:
(i) the baseline Q-GMRES method without preconditioning, and
(ii) a preconditioned variant using the built-in \emph{left LU preconditioner}.
This preconditioner is based on a quaternion LU factorization with partial pivoting and is one of the
practical solver enhancements contributed in \texttt{QuatIca}, designed to accelerate Krylov convergence
while preserving the native quaternion formulation.

The benchmark script generates several families of test matrices of increasing size, including:
\begin{itemize}[leftmargin=1.2em]
\item \emph{SPD}: Hermitian positive definite systems with controlled spectrum;
\item \emph{DENSE}: dense random quaternion systems with a stabilizing diagonal shift;
\item \emph{ILL-COND}: moderately ill-conditioned Hermitian systems with prescribed spectral decay.
\end{itemize}
For each size, scenario, and random seed, the script solves the corresponding system with both methods
and records:
\begin{itemize}[leftmargin=1.2em]
\item the number of iterations to convergence,
\item the wall-clock solve time,
\item the relative solution error
\[
\frac{\|\hat x-x_\star\|_F}{\|x_\star\|_F},
\]
\item the final residual reported by the solver,
\item and a success flag indicating whether the prescribed tolerance was reached within the iteration budget.
\end{itemize}

\paragraph{Demonstration script.}
The benchmark can be launched directly from the repository root via the analysis wrapper:
\begin{lstlisting}[language=bash,caption={Q-GMRES benchmark with and without LU preconditioning.},label={lst:qgmres_bench}]
# from the repo root
python run_analysis.py qgmres_bench
\end{lstlisting}

The script produces a full publication-style dashboard
(\texttt{qgmres\_final\_performance\_report.png}) together with individual PDF panels exported from the dashboard.

\begin{figure}[ht!]
  \centering

  \begin{subfigure}{0.48\linewidth}
    \centering
    \includegraphics[width=\linewidth]{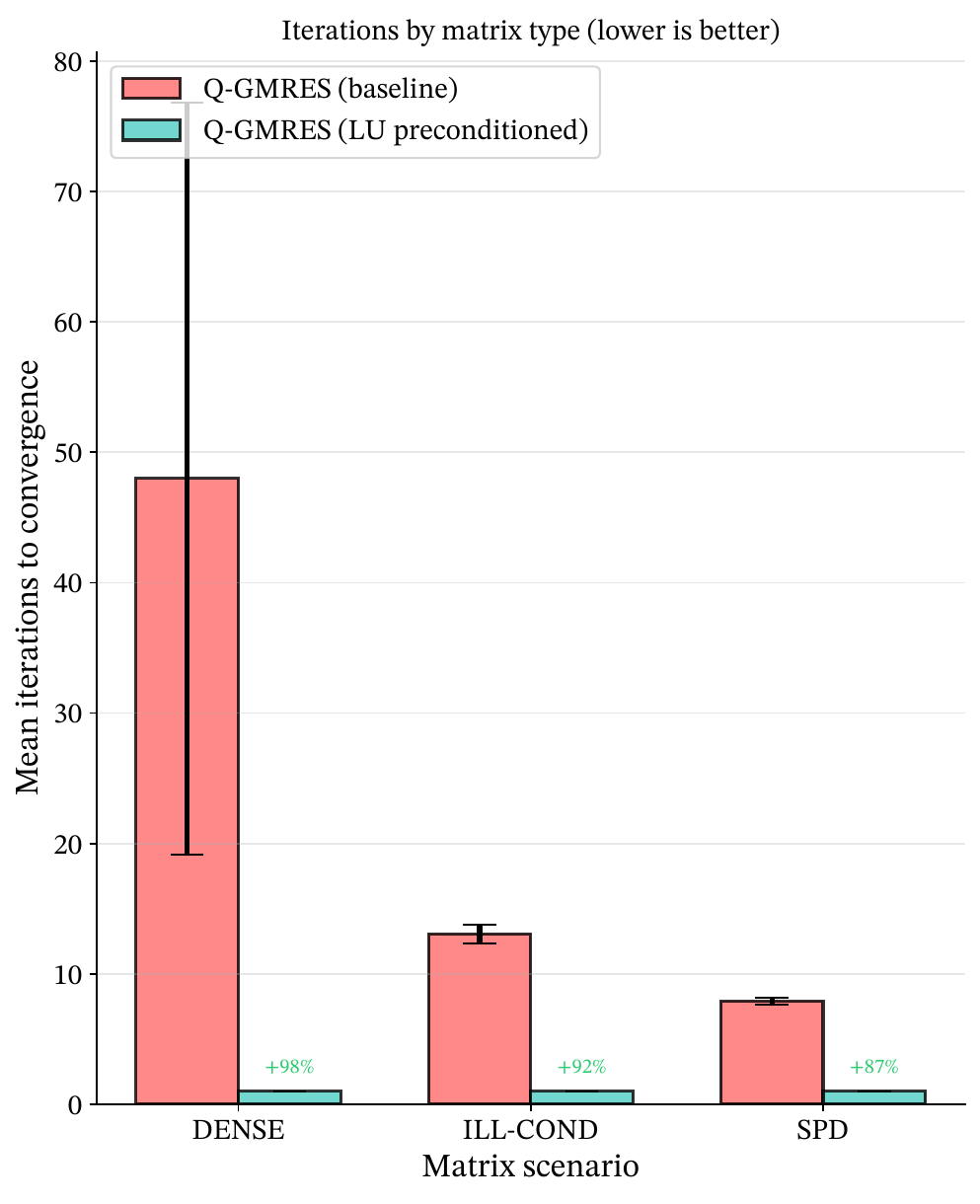}
    \caption{Mean number of iterations to convergence across matrix scenarios.}
    \label{fig:qgmres_iterations_by_scenario}
  \end{subfigure}\hfill
  \begin{subfigure}{0.48\linewidth}
    \centering
    
    \includegraphics[width=\linewidth]{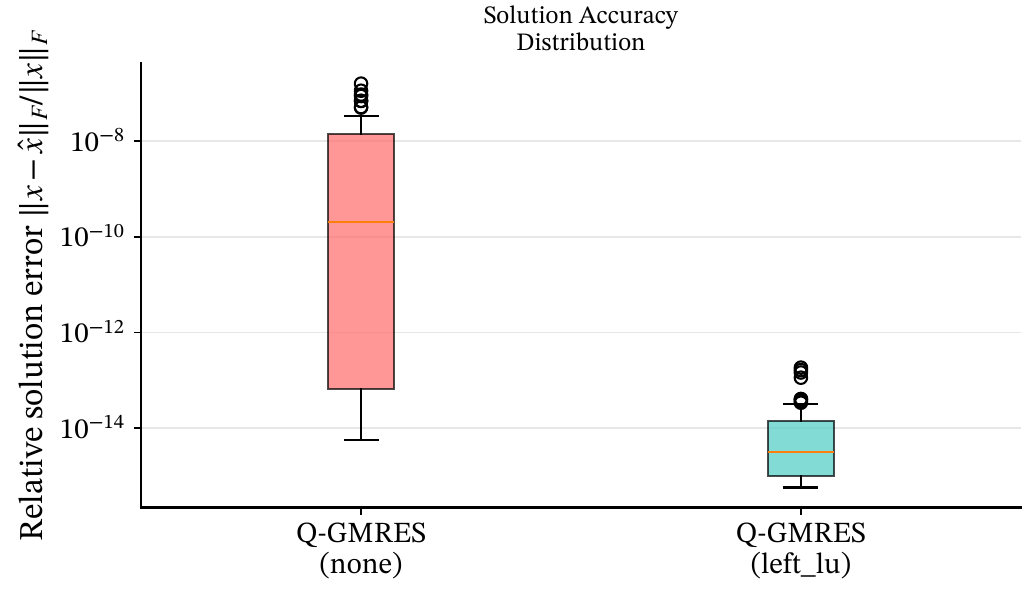}
    \caption{Distribution of relative solution errors for the two solver variants.}
    \label{fig:qgmres_accuracy_boxplot}
  \end{subfigure}

  \vspace{0.6em}

  \begin{subfigure}{0.90\linewidth}
    \centering
    \includegraphics[width=\linewidth]{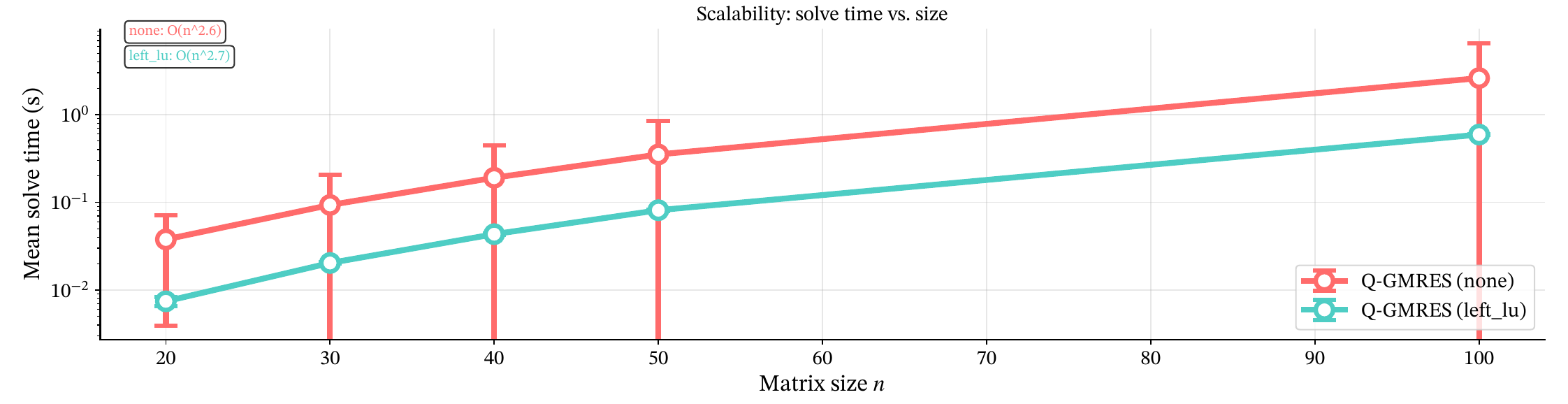}
    \caption{Solve-time scalability as the matrix dimension increases.}
    \label{fig:qgmres_scalability}
  \end{subfigure}

  \caption{Representative Q-GMRES benchmark outputs. The built-in LU preconditioner reduces iteration counts
  and significantly improves runtime, while both preconditioned and unpreconditioned Q-GMRES retain high solution accuracy, measured by the relative error with respect to \(x_\star\).}
  \label{fig:qgmres_benchmark}
\end{figure}

\paragraph{What the benchmark shows.}
The results highlight the practical value of the built-in LU preconditioner.
Across the tested scenarios, preconditioning consistently reduces the number of GMRES iterations required for
convergence and leads to noticeably lower solve times, with the gains becoming more visible as the matrix size
and difficulty of the problem increase. A representative subset of the benchmark outputs is shown in
Figure~\ref{fig:qgmres_benchmark}: Figures~\ref{fig:qgmres_iterations_by_scenario}
and~\ref{fig:qgmres_scalability} illustrate the reduction in iteration count and the improved runtime scaling,
whereas Figure~\ref{fig:qgmres_accuracy_boxplot} shows that both solver variants retain high solution accuracy,
as measured by the relative error with respect to the known reference solution \(x_\star\).
Overall, these results indicate that quaternion LU preconditioning is an effective and practically useful
accelerator for Q-GMRES within \texttt{QuatIca}.

\section{Current limitations and future directions}
\label{sec:outlook}
While \texttt{QuatIca} already covers a broad set of quaternion numerical routines, several
important directions remain open, both in terms of scalability and algorithmic coverage.

\paragraph{Scalability of conic optimization (\texttt{OptiQ}).}
Current limitations primarily concern scalability for large SDPs in \texttt{OptiQ}, where Newton steps assemble and factor dense Schur complements. Near-term improvements include:
(i) matrix-free Schur operators and iterative linear solves (CG/MINRES in $\R^m$ with quaternion-valued
operator applications),
(ii) caching and reuse of factorizations/spectral decompositions across line-search and residual evaluations,
(iii) robust Slater initialization strategies (projection--shift--snap and Phase-I feasibility formulations),
(iv) expanded ADMM variants tailored to large-scale and structured settings, and
(v) chordal decomposition techniques for sparse SDP instances. In the latter approach, when the aggregate sparsity pattern of the Hermitian matrix variable is sparse and admits a chordal extension, a large positive semidefinite constraint can be reformulated as a collection of smaller semidefinite constraints on clique submatrices together with overlap-consistency constraints. This can significantly reduce the cost of interior-point
steps and is therefore a particularly promising direction for large sparse problems.

\paragraph{Sparse and structured quaternion formats.}
Many target applications (imaging, graphs, inverse problems) exhibit structure and sparsity that are not yet systematically exploited. Extending the library with sparse quaternion operators and structured linear algebra (Toeplitz/circulant, convolutional, block-Hankel, and Kronecker forms) will enable larger problem sizes and improve performance across solvers.

\paragraph{Maximum-volume selection (quaternion MaxVol) and downstream applications.}
Efficient and scalable maximum-volume submatrix selection in the quaternion setting is a key missing building block
for robust low-rank approximation, cross/skeleton decompositions, CUR-type factorizations, and adaptive sampling.
Future work includes fast quaternion MaxVol and rect-MaxVol variants, careful treatment of noncommutativity, and
applications to data selection, compression, and preconditioning.

\paragraph{Randomized iterative solvers for quaternion linear systems.}
Beyond deterministic Krylov solvers, randomized techniques can reduce cost per iteration and improve robustness in
ill-conditioned regimes. Promising directions include randomized orthogonalization in Arnoldi/GMRES, sketching-based
preconditioners, and pass-efficient randomized residual estimation, leading to practical randomized quaternion GMRES
variants and hybrid deterministic-randomized solvers.

\paragraph{Quaternion tensor decompositions.}
Many modern datasets are naturally tensor-valued (video, dynamic graphs, multi-detector streams) and benefit from
compressed representations. Extending \texttt{QuatIca} with quaternion tensor decompositions, including Tucker~\cite{ya2025parallelizablequaternionhigherordersingular,MIAO2020107505},
CPD, tensor-train (TT)~\cite{MIAO2024111222} and tensor-ring~\cite{Miao2024} models, would provide a unified framework for compression and learning in the quaternion
domain. This also aligns with emerging applications such as low-rank Hankel tensor denoising and multi-channel
spatio-temporal inverse problems.

\paragraph{Advanced signal processing.}
We aim to extend the library with quaternion-valued signal processing tools for multicomponent data (e.g., RGB images, vector-sensor signals, dual-polarized communications). A first target is the \emph{Quaternion Fourier Transform} (QFT) and related time-frequency transforms (e.g., quaternion STFT / localized transforms), which allow holistic processing of coupled channels rather than channel-wise pipelines \cite{EllLeBihanSangwine2014QFT,AssefaEtAl2010LocalQFT}. A second target is quaternion adaptive filtering and widely-linear modeling, including QLMS and its widely-linear extensions, properness/impropriety diagnostics, and fast frequency-domain implementations, which are key ingredients for estimation and tracking with quaternion random signals \cite{TookMandic2009QLMS,TookMandic2010QWL,ViaRamirezSantamaria2010Properness,OrtolaniEtAl2017FDQAF}. In telecommunications, we will focus on dual-polarized MIMO channel modeling and equalization using quaternion arithmetic, including lattice-reduction-aided equalizers and blind widely-linear quaternion Godard/CMA-type algorithms \cite{SternFischer2018ICTDualPol,ZhangEtAl2025WLQGodard,AliEtAl2021QuaternionCodes}. Finally, we plan to include sparsity-driven recovery methods for quaternion signals (compressed sensing / sparse representation), such as quaternion \(\ell_1\)-minimization and reweighted sparse learning, to support reconstruction from limited (and potentially noisy) measurements \cite{GomesHartmannKahler2017CSQuaternion,WuEtAl2012L1Quaternion,ZouEtAl2023AWQSR}.

\paragraph{Longer-term perspective.}
A broader roadmap includes improved GPU acceleration, deeper integration with the scientific Python ecosystem, and
a wider range of reproducible application pipelines, such as polarimetric imaging, multi-sensor fusion, and
physics-informed quaternion models. Looking further ahead, the reference \texttt{MATLAB} QTFM toolbox also
supports octonion matrices in addition to quaternion ones. While \texttt{QuatIca} currently focuses on quaternion
numerical linear algebra, extending the library to octonion matrices represents a natural longer-term direction.

\section*{Acknowledgments}
\texttt{QuatIca} was inspired by pioneering work on quaternion numerical linear algebra, in particular the
QTFM (Quaternion Toolbox for MATLAB) developed by Stephen J.\ Sangwine and Nicolas Le Bihan. Their
comprehensive MATLAB implementation helped demonstrate the practical value of quaternion-based numerical
methods and provided a valuable reference point for several design and implementation choices in
\texttt{QuatIca}. We also gratefully acknowledge Roland Hildebrand for discussions that motivated the
development of \texttt{OptiQ} and, in particular, for proposing the idea of extending conic programming tools
to quaternion Hermitian optimization problems.

\bibliographystyle{plain}
\bibliography{quatica_refs}


\appendix

\section{OptiQ Newton system and Schur complement derivation}
\label{app:optiq_details}

This appendix collects the technical details underlying Algorithm~\ref{alg:optiq_barrier_path}.

\subsection{Barrier KKT system}
For fixed $\mu>0$, the barrier problem \eqref{eq:optiq_barrier_subproblem} has KKT conditions
\[
A(X)=b,\qquad C + A^*(y) - \mu X^{-1}=0,\qquad X\succ 0.
\]
Define residuals
\[
r_p := b-A(X),\qquad r_d := C + A^*(y) - \mu X^{-1}.
\]

\subsection{Newton linearization}
Let $H$ denote the Hessian operator of $-\mu\log\det(X)$:
\[
H[W] = \mu X^{-1} W X^{-1},
\qquad
H^{-1}[G] = \frac{1}{\mu}XGX.
\]
Newton's equations for $(\Delta X,\Delta y)$ are
\[
H[\Delta X] + A^*(\Delta y) = -r_d,\qquad A(\Delta X)=r_p.
\]
Eliminating $\Delta X$ yields the Schur complement system
\begin{equation}
\label{eq:optiq_schur}
(AH^{-1}A^*)\,\Delta y = -\,r_p - A\!\left(H^{-1}[r_d]\right),
\end{equation}
and then
\begin{equation}
\label{eq:optiq_backsub}
\Delta X = -\,H^{-1}\!\left[r_d + A^*(\Delta y)\right].
\end{equation}

\subsection{Hat-space form}
Using the hat operators \eqref{eq:optiq_hat_ops}, the Schur system \eqref{eq:optiq_schur} becomes
\[
(\hat A H^{-1}\hat A^*)\,\Delta \hat y = -\,\hat r_p - \hat A\!\left(H^{-1}[\hat r_d]\right),
\]
with $\hat r_p=\hat b-\hat A(X)$ and $\hat r_d=C+\hat A^*(\hat y)-\mu X^{-1}$, and the back-substitution is
$\Delta X = -H^{-1}[\hat r_d+\hat A^*(\Delta \hat y)]$.
OptiQ assembles the Schur matrix columnwise in $\R^{m\times m}$ and solves it with dense linear algebra
(with a tiny adaptive diagonal jitter if needed).

\subsection{Feasibility correction and line-search}
OptiQ uses:
(i) a fraction-to-boundary step rule to preserve positive definiteness of $X+t\Delta X$,
(ii) Armijo backtracking on a combined residual norm, and
(iii) the equality correction \eqref{eq:optiq_equality_snap} to suppress numerical drift in feasibility.

\section{Quaternion Cholesky factorization: dense and sparse implementations}
\label{app:chol_details}

For a quaternion Hermitian positive definite matrix \(A\in\Herm_n(\HH)\), the Cholesky factorization takes the form
\[
A = L L^H,
\]
where \(L\in\HH^{n\times n}\) is lower triangular with strictly positive real diagonal. In \texttt{QuatIca}, dense and sparse inputs are handled differently.

\paragraph{Dense native quaternion factorization.}
The routine \texttt{chol\_quat\_dense} implements a native quaternion analogue of the classical Cholesky algorithm. For a quaternion Hermitian positive definite matrix \(A\in\Herm_n(\HH)\), with entries \(a_{ij}\), the \(k\)-th pivot is computed as
\[
s_k = a_{kk} - \sum_{j<k} L_{kj}\,\overline{L_{kj}}
      = a_{kk} - \sum_{j<k} |L_{kj}|^2 \in \mathbb{R},
\]
and the diagonal entry is
\[
L_{kk} = \sqrt{s_k} \in \mathbb{R}_{>0}.
\]
For \(i>k\), the subdiagonal entries are updated through
\[
t_{ik} = a_{ik} - \sum_{j<k} L_{ij}\,\overline{L_{kj}},
\qquad
L_{ik} = \frac{t_{ik}}{L_{kk}}.
\]
Since \(L_{kk}\) is real and positive, this division is unambiguous and numerically safe. The implementation checks that the diagonal pivots are approximately real, optionally symmetrizes the input matrix, and allows a small diagonal \texttt{jitter} for stabilization.

Given the factor \(L\), the companion routine \texttt{solve\_chol\_quat\_dense} solves
\[
Ax=b,\qquad A=LL^H,
\]
by first solving \(Ly=b\) and then \(L^H x=y\). The implementation supports both vector and matrix right-hand sides.

\paragraph{Sparse implementation via complex embedding.}
For sparse matrices, \texttt{QuatIca} avoids implementing a full sparse quaternion triangular factorization directly. Instead, it uses the standard complex adjoint (or symplectic) embedding
\[
\chi(A)=
\begin{bmatrix}
X & Y\\
-\overline{Y} & \overline{X}
\end{bmatrix},
\qquad
X=A_w+iA_x,\quad Y=A_y+iA_z,
\]
which maps a quaternion matrix \(A=A_w + A_x \mathtt{i} + A_y \mathtt{j} + A_z \mathtt{k}\) to a complex matrix of dimension \(2n\times 2n\). For quaternion Hermitian positive definite matrices, this embedding preserves Hermiticity and positive definiteness. The sparse routine \texttt{chol\_quat\_sparse} factors \(\chi(A)\) with CHOLMOD through \texttt{scikit-sparse}, optionally after adding a small diagonal \texttt{jitter}.

Rather than returning an explicit quaternion triangular factor, the routine returns a factor object that solves quaternion linear systems through the embedded complex factorization and, when supported by the backend, evaluates a logarithmic determinant. In particular, for quaternion Hermitian positive definite matrices one has
\[
\log\det(A)=\frac12\log\det(\chi(A)),
\]
which is the identity used internally by the \texttt{logdet()} method of the sparse factor object.


\end{document}